\documentclass[11pt,reqno]{amsart}
\usepackage{amsmath}
\usepackage{amsmath,amssymb,amsfonts,amsthm}
\usepackage{hyperref}
\usepackage{mathrsfs}
\usepackage{enumitem}
\usepackage{geometry}
\hypersetup{
	colorlinks=true,
	linkcolor=blue,
	citecolor=blue,
	urlcolor=blue
}

\newtheorem{theorem}{Theorem}[section]

\newtheorem{proposition}[theorem]{Proposition}
\newtheorem{corollary}[theorem]{Corollary}

\theoremstyle{definition}
\newtheorem{definition}[theorem]{Definition}

\theoremstyle{remark}
\newtheorem{remark}[theorem]{Remark}

\title{An \(m\)-Hessian Approach to the Strominger-Yau-Zaslow Conjecture for Degenerating Calabi-Yau Manifolds}

\author{Truong Dinh Dat}

\address{Nguyen Trai University, Hanoi, Viet Nam}

\email{truongdinhdat14081994qb@gmail.com}

\begin{document}
	\begin{abstract}
	We investigate a pluripotential-theoretic approach to the Strominger--Yau--Zaslow (SYZ) conjecture for degenerating families of compact Calabi--Yau manifolds. Our approach is based on complex \(m\)-Hessian equations and their associated capacity and stability theories. The main objective is to explore how uniform pluripotential estimates may contribute to the analytical study of Ricci-flat metric degeneration and the geometric structures expected in the SYZ picture.
	
	For a family of Calabi--Yau manifolds approaching a maximally degenerate limit, we consider uniform \(m\)-Hessian capacity estimates on regions corresponding to the regular part of the essential skeleton. We discuss their potential applications to stability estimates, energy bounds, and compactness properties of complex Hessian potentials. We then examine the relationship between these analytical estimates and the Ricci-flat complex Monge--Ampère equations governing the Kähler metrics on the smooth fibers. Under suitable geometric and analytic assumptions, we formulate a framework for studying the limiting real Monge--Ampère geometry and the construction of local special Lagrangian torus fibrations.
	
	We further discuss the discriminant locus, monodromy, and the dualization of regular special Lagrangian torus fibers. The resulting semi-flat mirror geometry is supplemented by the holomorphic-disk contributions and wall-crossing transformations required for a quantum-corrected mirror construction. The role of complex \(m\)-Hessian theory is auxiliary: it provides a potential source of quantitative estimates, while the construction of special Lagrangian fibrations and global mirror spaces requires additional differential-geometric and symplectic-geometric input.
	\end{abstract}

	

	\maketitle	
\section{Introduction}

The Strominger--Yau--Zaslow (SYZ) conjecture is one of the central
geometric proposals underlying mirror symmetry. In its original
formulation, the conjecture suggests that a pair of mirror
Calabi--Yau manifolds should admit dual special Lagrangian torus
fibrations over a common base. More precisely, if
\[
f:X\longrightarrow B,
\qquad
\check f:\check X\longrightarrow B
\]
are special Lagrangian torus fibrations, then, over the regular
locus of the base, the fibers of the mirror fibration should be
obtained by dualizing the corresponding real tori. This geometric
picture was proposed by Strominger, Yau, and Zaslow
\cite{SYZ1996}.

The special Lagrangian condition is closely related to calibrated
geometry. A submanifold \(L\) of a Calabi--Yau manifold
\((X,\omega,\Omega)\) is special Lagrangian of phase
\(\vartheta\) if
\[
\omega|_L=0,
\qquad
\operatorname{Im}(e^{-i\vartheta}\Omega)|_L=0,
\qquad
\operatorname{Re}(e^{-i\vartheta}\Omega)|_L
=\operatorname{vol}_L.
\]
The general theory of calibrated submanifolds provides the
geometric framework for this condition
\cite{HarveyLawson1982}. Moreover, McLean's deformation theory
shows that infinitesimal deformations of a compact special
Lagrangian submanifold are governed by harmonic one-forms
\cite{McLean1998}. These results explain why the moduli spaces of
special Lagrangian submanifolds naturally carry affine structures
and why the duality of torus fibers is a plausible geometric
mechanism for mirror symmetry. Hitchin's work on the moduli space
of special Lagrangian submanifolds further clarifies this
relationship \cite{Hitchin1997}.

A major difficulty in the SYZ conjecture is the appearance of
singular fibers. In a degenerating family of Calabi--Yau manifolds,
the Ricci-flat metrics are expected to collapse, after suitable
rescaling or normalization, to a lower-dimensional space equipped
with an affine structure and a metric satisfying a real
Monge--Amp\`ere equation. The Gross--Wilson analysis of degenerating
K3 surfaces provides a fundamental model for this phenomenon
\cite{GrossWilson2000}. Joyce's work emphasizes that singularities
of special Lagrangian fibrations are subtle and that the
discriminant locus may have a complicated structure
\cite{Joyce2003}. Consequently, constructing a global special
Lagrangian fibration requires more than the local deformation
theory of individual fibers.

A complementary direction was developed by Kontsevich and Soibelman \cite{KontsevichSoibelman2006}, who proposed a non-Archimedean approach to mirror symmetry. Their construction relates integral affine manifolds with singularities to analytic spaces over non-Archimedean fields. In this framework, the limiting affine geometry of a degenerating Calabi--Yau manifold is connected with its non-Archimedean analytification. The essential skeleton and its associated integral affine structures provide a natural candidate for the base of the limiting torus fibration. This viewpoint also links the metric collapse of Calabi--Yau manifolds with the algebraic and non-Archimedean geometry of degenerations.

Gross and Siebert developed a systematic algebro-geometric program for mirror symmetry using logarithmic geometry and toric degenerations \cite{GrossSiebert2006}. Their approach associates affine manifolds with singularities and polyhedral decompositions to suitable degenerations of Calabi--Yau varieties. A discrete Legendre transform then produces the combinatorial data of the mirror degeneration. The subsequent development of scattering diagrams and wall structures provides a mechanism for correcting the naive gluing of semi-flat mirror charts. This program has become one of the principal mathematical approaches to constructing mirrors from degeneration data, particularly in settings where a direct construction of special Lagrangian fibrations is difficult. In particular,
weight functions on non-Archimedean analytic spaces provide a
systematic construction of the relevant skeleton
\cite{MustataNicaise2015}, while the essential skeleton gives a
birationally meaningful object associated with the degeneration
\cite{NicaiseXu2016}. These constructions are closely related to
the expected base of the SYZ fibration, although identifying the
essential skeleton with a metric limit requires additional
geometric and analytic arguments.

Another essential aspect of SYZ mirror symmetry is the appearance of quantum corrections. The naive dualization of torus fibers generally produces only a semi-flat mirror, which fails to capture the full geometry near the discriminant locus. Holomorphic disks bounded by Lagrangian torus fibers contribute additional terms to the complex structure and the gluing transformations of the mirror. Auroux \cite{Auroux2007,Auroux2009} developed important examples and conceptual explanations of wall-crossing phenomena through Floer theory and disk counting. In this perspective, the discontinuity of disk counts across walls is compensated by nontrivial coordinate transformations, leading to corrected mirror spaces. These ideas are closely related to the scattering-diagram formalism of Gross and Siebert.

The analytical study of degenerating Calabi--Yau metrics is
closely connected with the complex Monge--Amp\`ere equation.
By Yau's solution of the Calabi conjecture, a compact K\"ahler
manifold with vanishing first Chern class admits a Ricci-flat
K\"ahler metric in each K\"ahler class \cite{Yau1978}. The
construction of weak solutions and the pluripotential theory of
complex Monge--Amp\`ere equations were developed through the
foundational work of Bedford and Taylor
\cite{BedfordTaylor1976} and the global estimates of
Ko\l odziej \cite{Kolodziej1998}. These theories form an essential
part of the analytical background for studying Ricci-flat
potentials and their possible degeneration.

In this paper, we investigate a complementary analytical
framework based on complex \(m\)-Hessian equations. For
\(1\leq m<n\), the complex \(m\)-Hessian operator is given
formally by
\[
H_m(u)
=
(dd^c u)^m\wedge\omega^{n-m},
\]
where \((X,\omega)\) is a compact K\"ahler manifold and \(u\)
belongs to an appropriate class of \(m\)-subharmonic functions.
The theory of weak solutions to complex Hessian equations was
developed by B{\l}ocki \cite{Blocki2005}. A priori estimates
and uniqueness results were subsequently established by
Dinew and Ko\l odziej \cite{DinewKolodziej2014} and by Dinew
and Lu \cite{DinewLu2015}. Further developments include weak
solutions on compact Hermitian manifolds
\cite{KolodziejNguyen2016}, Liouville and Calabi--Yau type
theorems for complex Hessian equations
\cite{DinewKolodziej2017}, and the study of non-standard
properties of \(m\)-subharmonic functions
\cite{DinewKolodziej2018}. In addition, Ko\l odziej and
Nguyen established an inequality relating complex Hessian
measures of H\"older continuous \(m\)-subharmonic functions
to Hessian capacity \cite{KolodziejNguyen2020}.

These developments suggest that Hessian capacity and
pluripotential stability may provide useful tools for studying
families of degenerating complex manifolds. In particular,
capacity estimates can quantify the size of exceptional sets,
while stability estimates can measure the dependence of
solutions on their Monge--Amp\`ere or Hessian data. Such estimates
may be useful in obtaining uniform control of auxiliary
potentials on suitable regions of a degenerating family.
However, the complex \(m\)-Hessian equation with \(m<n\) is not
the equation defining the Ricci-flat metric: the latter is
governed by the complex Monge--Amp\`ere equation of order \(n\).
Thus, the role of the \(m\)-Hessian theory in the present work
is complementary. It is intended to provide additional
pluripotential estimates that may assist the analysis of
degenerating Ricci-flat metrics, rather than to replace the
Calabi--Yau equation.

The metric and analytic aspects of the SYZ conjecture have also continued to develop. In particular, Yang Li \cite{Li2022} established a major result for Calabi--Yau hypersurfaces in the Fermat family, proving the expected metric collapse in that setting and obtaining a detailed description of the limiting geometry. His work provides an important higher-dimensional contribution to the analytic SYZ program and illustrates the extent to which the conjectural picture can be established in geometrically structured families. Further work has developed the relationship
between the metric SYZ conjecture and non-Archimedean geometry
\cite{Li2023}, as well as detailed local geometric models
involving Taub--NUT and Ooguri--Vafa type metrics
\cite{Li2023Memoirs}. Nevertheless, the general problem of constructing and understanding special Lagrangian torus fibrations, including their singular fibers and quantum-corrected mirrors, remains open in broad classes of degenerating Calabi--Yau manifolds of complex dimension \(n\geq 3\).

Motivated by these developments, we propose an
\(m\)-Hessian-based pluripotential framework for investigating
degenerating Calabi--Yau manifolds and the analytical problems
associated with the SYZ picture. We consider a family
\[
\pi:\mathcal X\longrightarrow\mathbb D
\]
whose smooth fibers \(X_t\), for \(t\neq0\), are compact
Calabi--Yau manifolds of complex dimension \(n\). Our approach
focuses on the interaction between Hessian capacities,
stability estimates for complex Hessian equations, and the
complex Monge--Amp\`ere equations defining the Ricci-flat
metrics on \(X_t\).

The main objectives of this framework are as follows.

\begin{enumerate}
	\item We formulate uniform \(m\)-Hessian capacity estimates
	for suitable subsets of the degenerating fibers and identify
	the geometric assumptions needed to keep the relevant
	constants independent of the degeneration parameter.
	
	\item We investigate quantitative stability estimates for
	families of complex \(m\)-Hessian equations with varying
	right-hand sides. The purpose is to obtain uniform control
	of the associated potentials under appropriate energy and
	capacity assumptions.
	
	\item We study how these pluripotential estimates may be
	combined with the complex Monge--Amp\`ere equations for
	Ricci-flat metrics and with additional geometric hypotheses
	governing metric collapse.
	
	\item On the regular part of a candidate limiting base, we
	examine the conditions under which the limiting metric
	admits a Hessian-potential description and the local
	geometry is compatible with special Lagrangian torus
	fibrations.
	
	\item We discuss the role of the discriminant locus,
	monodromy, and dual torus fibrations in the construction
	of a candidate mirror space, together with the additional
	quantum corrections required for a global mirror
	construction.
\end{enumerate}

We emphasize that the purpose of this paper is to develop an
analytical framework and a collection of conditional geometric
steps toward the SYZ conjecture. Uniform \(m\)-Hessian
estimates alone do not imply metric collapse, the existence of
special Lagrangian fibrations, or the construction of a
quantum-corrected mirror. Each of these conclusions requires
additional geometric, analytic, and topological arguments.
The central question is whether uniform pluripotential
estimates can provide a useful bridge between the analysis of
degenerating complex Hessian equations and the geometric
structures expected in the SYZ picture.


\section{Preliminaries}
\label{sec:preliminaries}

In this section, we introduce the basic notation and geometric
structures used throughout the paper. We first recall the relevant
notions from Kähler geometry and complex Hessian potential theory.
We then describe the degeneration setting, the associated
non-Archimedean geometry, and the geometric structures underlying
the Strominger--Yau--Zaslow conjecture.

\subsection{Kähler geometry and complex differential operators}

Let \(X\) be a compact connected complex manifold of complex
dimension \(n\). We use the convention
\[
dd^c=\sqrt{-1}\,\partial\bar\partial.
\]
Let \(\theta\) be a Kähler form on \(X\). For a real-valued
function \(u\in C^2(X,\mathbb R)\), we write
\[
\theta_u=\theta+dd^c u.
\]
The function \(u\) is called \(\theta\)-plurisubharmonic if
\[
\theta_u\geq 0
\]
in the sense of real \((1,1)\)-forms. The space of such functions
is denoted by
\[
\operatorname{PSH}(X,\theta)
=
\left\{
u\in L^1(X):
u\text{ is upper semicontinuous and }
\theta+dd^c u\geq0
\right\}.
\]

For a smooth real \((1,1)\)-form \(\alpha\), let
\[
\lambda(\alpha)
=
(\lambda_1,\ldots,\lambda_n)
\]
denote its eigenvalues with respect to \(\theta\), so that
\[
\alpha\wedge\theta^{n-1}
=
\frac{\sigma_1(\lambda)}{n}\theta^n,
\]
where \(\sigma_k\) denotes the \(k\)-th elementary symmetric
polynomial. More generally,
\[
\alpha^k\wedge\theta^{n-k}
=
\frac{\sigma_k(\lambda)}{\binom nk}\theta^n.
\]

For \(1\leq m\leq n\), define the Gårding cone
\[
\Gamma_m
=
\left\{
\lambda\in\mathbb R^n:
\sigma_k(\lambda)\geq0,\quad
1\leq k\leq m
\right\}.
\]
A smooth real \((1,1)\)-form \(\alpha\) is called
\(m\)-positive with respect to \(\theta\) if
\[
\lambda(\alpha)\in\overline{\Gamma_m}.
\]

A function \(u\in C^2(X,\mathbb R)\) is called
\(\theta\)-\(m\)-subharmonic if
\[
\theta+dd^c u
\]
is \(m\)-positive with respect to \(\theta\). We denote the
corresponding class of upper semicontinuous \(m\)-subharmonic
functions by
\[
\operatorname{SH}_m(X,\theta).
\]
In particular,
\[
\operatorname{SH}_n(X,\theta)
=
\operatorname{PSH}(X,\theta).
\]

For sufficiently regular \(u\), the complex \(m\)-Hessian
operator is
\[
H_m(u)
=
(\theta+dd^c u)^m\wedge\theta^{n-m}.
\]
When \(u\) is nonsmooth, this operator is understood in the
pluripotential sense whenever the corresponding Hessian measure
is well-defined.

\subsection{Complex \(m\)-Hessian equations}

Let \(1\leq m\leq n\), and let
\(f\geq0\) be a measurable function on \(X\). The complex
\(m\)-Hessian equation takes the form
\begin{equation}
	\label{eq:complex-m-hessian}
	(\theta+dd^c u)^m\wedge\theta^{n-m}
	=
	f\,\theta^n.
\end{equation}
For a smooth strictly \(m\)-subharmonic solution, this is a
fully nonlinear elliptic equation in the admissible cone
\(\Gamma_m\). The case \(m=1\) is linear, up to normalization,
whereas the case \(m=n\) is the complex Monge--Ampère equation.

On a compact Kähler manifold, integration of
\eqref{eq:complex-m-hessian} gives the necessary compatibility
condition
\[
\int_X f\,\theta^n
=
\int_X\theta^n,
\]
provided \(u\) is sufficiently regular. Indeed, by Stokes'
theorem,
\[
\int_X(\theta+dd^c u)^m\wedge\theta^{n-m}
=
\int_X\theta^n.
\]

Since the equation is invariant under addition of constants to
\(u\), we impose the normalization
\[
\sup_X u=0.
\]

In the setting of this paper, we are particularly interested
in the range
\[
1\leq m<n.
\]
The corresponding Hessian equations provide an auxiliary
hierarchy of nonlinear equations whose capacity and stability
properties will be used to investigate uniform pluripotential
estimates in degenerating families.

\subsection{Hessian capacity and energy classes}

Let \(E\subset X\) be a Borel set. The global
\(m\)-Hessian capacity associated with \(\theta\) is defined by
\begin{equation}
	\label{eq:m-hessian-capacity}
	\operatorname{Cap}_{m,\theta}(E)
	=
	\sup\left\{
	\int_E
	(\theta+dd^c v)^m\wedge\theta^{n-m}:
	v\in\operatorname{SH}_m(X,\theta),
	\ -1\leq v\leq0
	\right\}.
\end{equation}
For nonsmooth test functions, the integral is interpreted
using the corresponding Hessian measure.

The capacity
\(\operatorname{Cap}_{m,\theta}\) is monotone with respect
to set inclusion. Namely, if \(E\subset F\), then
\[
\operatorname{Cap}_{m,\theta}(E)
\leq
\operatorname{Cap}_{m,\theta}(F).
\]
Moreover,
\[
\operatorname{Cap}_{m,\theta}(\varnothing)=0,
\qquad
\operatorname{Cap}_{m,\theta}(X)<+\infty
\]
under the usual compact Kähler assumptions.

The capacity provides a natural way to quantify the size of
exceptional sets in complex Hessian potential theory. In
particular, it is useful for studying convergence in capacity,
quasi-continuity, and stability of weak solutions.

For \(p>0\), an associated finite-energy class can be described
formally by
\[
\mathcal E_m^p(X,\theta)
=
\left\{
u\in\operatorname{SH}_m(X,\theta):
\sup_Xu=0,\ 
\int_X(-u)^p H_m(u)<+\infty
\right\},
\]
where the precise definition for unbounded functions is
understood through the standard truncation procedure and the
appropriate extension of the Hessian operator.

For a family of normalized solutions \(u_t\), a typical
uniform energy condition is
\begin{equation}
	\label{eq:uniform-energy}
	\sup_{t}
	\int_{X_t}(-u_t)^p
	(\theta_t+dd^c u_t)^m\wedge\theta_t^{n-m}
	<+\infty.
\end{equation}
Such a bound is an additional analytical requirement in a
degenerating family; it does not follow solely from
\(L^1\)-convergence of the right-hand sides.

\subsection{Degenerating families of Calabi--Yau manifolds}

Let \(\mathbb D\subset\mathbb C\) be the unit disk and consider
a proper flat holomorphic map
\[
\pi:\mathcal X\longrightarrow\mathbb D.
\]
We assume that the restriction
\[
\pi:\mathcal X^*
\longrightarrow\mathbb D^*,
\qquad
\mathcal X^*=\pi^{-1}(\mathbb D^*),
\]
is a smooth family, where
\(\mathbb D^*=\mathbb D\setminus\{0\}\).

For \(t\in\mathbb D^*\), write
\[
X_t=\pi^{-1}(t).
\]
We assume that each \(X_t\) is a compact Kähler manifold
of complex dimension \(n\) with trivial canonical bundle.
Thus, each fiber admits a nowhere-vanishing holomorphic
\(n\)-form
\[
\Omega_t\in H^0(X_t,K_{X_t}).
\]

After a finite base change and a suitable semistable
modification, when necessary, we may work with a model whose
central fiber is a reduced simple normal crossing divisor
\[
\mathcal X_0
=
\sum_{i\in I}D_i.
\]
Locally near a point of the central fiber, the degeneration
can be expressed in the form
\[
t=z_0z_1\cdots z_r,
\qquad 1\leq r\leq n,
\]
in suitable holomorphic coordinates on the total space.

The dual intersection complex of the central fiber is the
polyhedral complex whose vertices correspond to the irreducible
components \(D_i\), and whose cells record their nonempty
intersections. The essential skeleton, denoted by
\[
B=\operatorname{Sk}_{\mathrm{ess}}(X_K),
\]
is a distinguished subcomplex of the dual intersection
complex, defined using the weight function associated with
pluricanonical forms on the non-Archimedean generic fiber
\(X_K\). Here \(K=\mathbb C((t))\), after fixing the
degeneration parameter.

We focus on the maximal degeneration situation
\begin{equation}
	\label{eq:maximal-degeneration}
	\dim_{\mathbb R}B=n.
\end{equation}
This condition identifies the essential skeleton as an
\(n\)-dimensional candidate for the base of the limiting
SYZ fibration.

We write
\[
B_0=B\setminus\Delta,
\]
where \(\Delta\subset B\) denotes the discriminant locus.
The regular locus \(B_0\) is expected, under suitable geometric
conditions, to carry an integral affine structure. The existence
of such a structure is part of the geometric information
required for the SYZ construction and is not a consequence of
the normal crossing condition alone.

\subsection{Ricci-flat metrics and the complex Monge--Ampère equation}

For each \(t\in\mathbb D^*\), let \(\theta_t\) be a Kähler
form on \(X_t\). Define the positive volume form
\[
dV_{\Omega_t}
=
(\sqrt{-1})^{n^2}
\Omega_t\wedge\overline{\Omega_t},
\]
with the orientation chosen so that \(dV_{\Omega_t}>0\).

Set
\[
c_t
=
\frac{\displaystyle\int_{X_t}\theta_t^n}
{\displaystyle\int_{X_t}dV_{\Omega_t}}.
\]
Yau's solution of the Calabi conjecture gives a unique
smooth function \(\varphi_t\), normalized by
\[
\sup_{X_t}\varphi_t=0,
\]
such that
\begin{equation}
	\label{eq:ricci-flat-ma}
	(\theta_t+dd^c\varphi_t)^n
	=
	c_t\,dV_{\Omega_t}.
\end{equation}
The form
\[
\omega_t=\theta_t+dd^c\varphi_t
\]
is a Ricci-flat Kähler form in the cohomology class
\([\theta_t]\).

Indeed, since \(\Omega_t\) is holomorphic and nowhere
vanishing, the volume form \(dV_{\Omega_t}\) has vanishing
Ricci form. Equation \eqref{eq:ricci-flat-ma} therefore
implies
\[
\operatorname{Ric}(\omega_t)=0.
\]

It is important to distinguish
\eqref{eq:ricci-flat-ma} from the complex \(m\)-Hessian
equation. The Ricci-flat metric is governed by the
Monge--Ampère operator of order \(n\), whereas the
\(m\)-Hessian operators with \(m<n\) will be used as
auxiliary analytical tools.

\subsection{Non-Archimedean geometry and the limiting base}

Let
\[
R=\mathbb C[[t]],
\qquad
K=\mathbb C((t)).
\]
The generic fiber of the degeneration determines a
non-Archimedean analytic space \(X_K^{\mathrm{an}}\).
A semistable model of \(X_K\) gives rise to a Berkovich
skeleton, which is a polyhedral subset of
\(X_K^{\mathrm{an}}\).

The essential skeleton
\[
B=\operatorname{Sk}_{\mathrm{ess}}(X_K)
\]
is defined intrinsically using the weight function of a
pluricanonical form. In the maximal degeneration case,
its dimension is \(n\), as in
\eqref{eq:maximal-degeneration}.

The non-Archimedean viewpoint provides a natural framework
for comparing the degeneration of complex fibers with the
geometry of the limiting base. In particular, measures
associated with holomorphic volume forms may admit limiting
descriptions on the essential skeleton. Such convergence,
when available, is a measure-theoretic statement and does
not by itself imply convergence of the Ricci-flat metrics
or the existence of a special Lagrangian fibration.

In the sequel, we use the regular part \(B_0\) as the
candidate base for the limiting torus fibration. Any
integral affine structure or metric structure on \(B_0\)
will be considered under the corresponding additional
regularity assumptions.

\subsection{Special Lagrangian submanifolds and dual torus fibers}

Let \((X,\omega,J,\Omega)\) be a Calabi--Yau manifold of
complex dimension \(n\), where \(\omega\) is a Kähler form
and \(\Omega\) is a nowhere-vanishing holomorphic
\(n\)-form compatible with the Ricci-flat metric.

A real \(n\)-dimensional submanifold \(L\subset X\) is
Lagrangian if
\[
\omega|_L=0.
\]
For a fixed phase \(\vartheta\in\mathbb R\), an oriented
Lagrangian submanifold \(L\) is special Lagrangian of
phase \(\vartheta\) if
\begin{equation}
	\label{eq:special-lagrangian}
	\operatorname{Im}
	\left(e^{-\sqrt{-1}\vartheta}\Omega\right)|_L=0,
	\qquad
	\operatorname{Re}
	\left(e^{-\sqrt{-1}\vartheta}\Omega\right)|_L>0.
\end{equation}
After a suitable normalization of \(\Omega\), such a
submanifold is calibrated by
\[
\operatorname{Re}
\left(e^{-\sqrt{-1}\vartheta}\Omega\right).
\]

The SYZ conjecture predicts that, near a suitable large
complex structure limit, a Calabi--Yau manifold admits a
special Lagrangian torus fibration over an \(n\)-dimensional
base. On the regular locus, this takes the form
\[
f:X^{\mathrm{reg}}\longrightarrow B_0,
\qquad
f^{-1}(b)\simeq T^n.
\]

Given a smooth Lagrangian torus fiber \(L_b\), its dual
torus is naturally defined by
\[
\check L_b
=
H^1(L_b,\mathbb R)/H^1(L_b,\mathbb Z).
\]
Equivalently,
\[
\check L_b
\simeq
\operatorname{Hom}
\left(H_1(L_b,\mathbb Z),U(1)\right).
\]
The union of these dual tori over the regular base gives
the semi-flat candidate mirror
\[
\check X_0
=
\bigcup_{b\in B_0}\check L_b.
\]
The extension of this space across the discriminant locus
typically requires additional corrections arising from
holomorphic disks and wall-crossing phenomena.

\subsection{Standing notation and conventions}

Throughout the paper, we use the following notation:
\begin{itemize}
	\item \(X_t\) denotes a smooth compact Calabi--Yau fiber
	of the degeneration.
	\item \(\theta_t\) denotes a reference Kähler form on
	\(X_t\).
	\item \(\omega_t=\theta_t+dd^c\varphi_t\) denotes the
	Ricci-flat Kähler form determined by
	\eqref{eq:ricci-flat-ma}.
	\item \(H_m(u)\) denotes the complex \(m\)-Hessian measure
	\[
	H_m(u)
	=
	(\theta_t+dd^c u)^m\wedge\theta_t^{n-m}.
	\]
	\item \(\operatorname{Cap}_{m,\theta_t}\) denotes the
	corresponding global \(m\)-Hessian capacity.
	\item \(B=\operatorname{Sk}_{\mathrm{ess}}(X_K)\) denotes
	the essential skeleton, and
	\(B_0=B\setminus\Delta\) its regular locus.
	\item \(\Delta\) denotes the discriminant locus of the
	expected limiting torus fibration.
\end{itemize}

All uniform estimates in the degeneration parameter are
understood to hold for \(0<|t|<t_0\), with constants
independent of \(t\), unless otherwise stated. When an
estimate is asserted only on compact subsets of \(B_0\),
this dependence will be explicitly indicated.

	\section{Construction of the Algebraic Degeneration Model}
	
	Let \(n\geq 3\), and let
	\[
	\pi:\mathcal{X}\longrightarrow \mathbb{D}
	\]
	be a proper flat holomorphic family over the unit disk
	\[
	\mathbb{D}=\{t\in\mathbb{C}:|t|<1\}.
	\]
	Assume that \(\mathcal{X}\) is a normal complex space and that
	the restriction
	\[
	\pi:\mathcal{X}^{*}\longrightarrow\mathbb{D}^{*},
	\qquad
	\mathcal{X}^{*}=\pi^{-1}(\mathbb{D}^{*}),
	\]
	is a smooth projective morphism, where
	\[
	\mathbb{D}^{*}=\mathbb{D}\setminus\{0\}.
	\]
	
	For every \(t\in\mathbb{D}^{*}\), denote
	\[
	X_t=\pi^{-1}(t).
	\]
	We assume that each fiber \(X_t\) is a compact Calabi--Yau
	manifold of complex dimension \(n\), namely,
	\[
	K_{X_t}\cong\mathcal{O}_{X_t},
	\qquad
	h^{i,0}(X_t)=0
	\quad\text{for }0<i<n.
	\]
	
	Moreover, assume that the family is polarized by a relatively
	ample line bundle \(\mathcal{L}\) on \(\mathcal{X}\).
	
	\subsection{1. Semistable Reduction}
	
	After a finite base change
	\[
	t=s^N,
	\qquad N\in\mathbb{N},
	\]
	and a suitable birational modification, we may work with a
	semistable model
	\[
	\pi:\mathcal{X}\longrightarrow\mathbb{D}
	\]
	whose central fiber is a reduced simple normal crossing
	divisor:
	\[
	\mathcal{X}_0
	=
	\pi^{-1}(0)
	=
	\sum_{i=1}^{N_0}D_i.
	\]
	
	Here each \(D_i\) is a smooth irreducible divisor, and locally
	near every point of \(\mathcal{X}_0\), there exist holomorphic
	coordinates
	\[
	(z_0,z_1,\ldots,z_n)
	\]
	such that
	\[
	\pi(z_0,\ldots,z_n)
	=
	z_0z_1\cdots z_r,
	\qquad 1\leq r\leq n+1.
	\]
	
	Thus the central fiber has simple normal crossings.
	
	We denote by
	\[
	D_I=\bigcap_{i\in I}D_i
	\]
	the intersection of the irreducible components indexed by
	\(I\subset\{1,\ldots,N_0\}\).
	
	The dual intersection complex
	\[
	\Delta(\mathcal{X}_0)
	\]
	is the simplicial complex whose vertices correspond to the
	components \(D_i\), and whose \(k\)-dimensional simplices
	correspond to nonempty intersections of \(k+1\) distinct
	components.
	
	\subsection{2. The Essential Skeleton}
	
	Let
	\[
	K=\mathbb{C}((t))
	\]
	be the field of formal Laurent series, and let
	\[
	X_K=\mathcal{X}\times_{\mathbb{D}}\operatorname{Spec}K
	\]
	be the generic fiber over \(K\).
	
	The Berkovich analytification of \(X_K\) is denoted by
	\[
	X_K^{\mathrm{an}}.
	\]
	
	The semistable model determines a skeleton
	\[
	\operatorname{Sk}(\mathcal{X})
	\subset X_K^{\mathrm{an}},
	\]
	which is naturally identified with a subcomplex of
	\[
	\Delta(\mathcal{X}_0).
	\]
	
	We define the essential skeleton by
	\[
	\operatorname{Sk}_{\mathrm{ess}}(X_K)
	\subseteq\operatorname{Sk}(\mathcal{X}).
	\]
	
	We impose the maximal degeneration condition
	\[
	\boxed{
		\dim_{\mathbb{R}}
		\operatorname{Sk}_{\mathrm{ess}}(X_K)=n.
	}
	\tag{1.1}
	\]
	
	This condition means that the essential skeleton has the
	maximal possible real dimension.
	
	In particular, it provides an \(n\)-dimensional candidate
	for the base space of the SYZ fibration.
	
	\subsection{3. Integral Affine Structure on the Regular Locus}
	
	Let
	\[
	B=\operatorname{Sk}_{\mathrm{ess}}(X_K)
	\]
	and let
	\[
	\Delta\subset B
	\]
	denote its discriminant locus.
	
	We set
	\[
	B_0=B\setminus\Delta.
	\]
	
	On the regular locus \(B_0\), we seek an integral affine
	structure defined by coordinate charts
	\[
	x_\alpha:U_\alpha\longrightarrow\mathbb{R}^n,
	\]
	whose transition functions belong to
	\[
	GL(n,\mathbb{Z})\ltimes\mathbb{R}^n.
	\]
	
	More precisely, for two overlapping charts, we require
	\[
	x_\beta
	=
	A_{\alpha\beta}x_\alpha+b_{\alpha\beta},
	\]
	where
	\[
	A_{\alpha\beta}\in GL(n,\mathbb{Z}),
	\qquad
	b_{\alpha\beta}\in\mathbb{R}^n.
	\]
	
	The existence and compatibility of this affine structure
	will be established under additional geometric hypotheses
	in the subsequent steps.
	
	\subsection{4. The Calabi--Yau Volume Form}
	
	Since the generic fibers are Calabi--Yau, we assume that
	the relative canonical bundle is trivial over the smooth
	locus:
	\[
	K_{\mathcal{X}/\mathbb{D}}
	\cong\mathcal{O}_{\mathcal{X}}.
	\]
	
	Thus there exists a nowhere-vanishing relative holomorphic
	\(n\)-form
	\[
	\Omega\in
	H^0(\mathcal{X}^{*},K_{\mathcal{X}/\mathbb{D}}).
	\]
	
	Its restriction to each fiber is denoted by
	\[
	\Omega_t=\Omega|_{X_t}.
	\]
	
	We normalize the corresponding volume measure by
	\[
	d\mu_t
	=
	\frac{(\sqrt{-1})^{n^2}
		\Omega_t\wedge\overline{\Omega_t}}
	{\displaystyle
		\int_{X_t}(\sqrt{-1})^{n^2}
		\Omega_t\wedge\overline{\Omega_t}}.
	\]
	
	The behavior of these measures as \(t\to0\) will be
	studied in the non-Archimedean setting.
	
	\subsection{5. Ricci-Flat Kähler Metrics}
	
	Let
	\[
	\theta_t\in c_1(\mathcal{L}|_{X_t})
	\]
	be a family of Kähler forms.
	
	We consider the rescaled Kähler classes
	\[
	[\omega_t]
	=
	\lambda_t[\theta_t],
	\qquad
	\lambda_t=|\log|t||.
	\]
	
	By Yau's solution of the Calabi conjecture, each class
	contains a unique Ricci-flat Kähler metric \(\omega_t\)
	after fixing its total volume.
	
	Consequently,
	\[
	\operatorname{Ric}(\omega_t)=0
	\]
	and
	\[
	\omega_t^n=c_t
	(\sqrt{-1})^{n^2}
	\Omega_t\wedge\overline{\Omega_t},
	\]
	where \(c_t>0\) is determined by the volume normalization.
	
	\subsection{6. Initial Geometric Setup for the SYZ Conjecture}
	
	Under the maximal degeneration condition (1.1), the
	essential skeleton
	\[
	B=\operatorname{Sk}_{\mathrm{ess}}(X_K)
	\]
	has real dimension \(n\).
	
	The SYZ conjecture predicts that, away from a suitable
	discriminant locus, the manifolds \(X_t\) admit
	special Lagrangian torus fibrations
	\[
	f_t:U_t\longrightarrow B_t,
	\]
	where
	\[
	U_t\subset X_t,
	\qquad
	B_t\subset B_0,
	\]
	and the generic fibers satisfy
	\[
	f_t^{-1}(b)\cong T^n.
	\]
	
	The special Lagrangian conditions are
	\[
	\omega_t|_{f_t^{-1}(b)}=0
	\]
	and
	\[
	\operatorname{Im}\Omega_t|_{f_t^{-1}(b)}=0.
	\]
	
	The purpose of the subsequent sections is to establish
	uniform pluripotential estimates, metric convergence,
	and the existence of these special Lagrangian torus
	fibrations.
	
	\begin{remark}
		The semistable model and the maximal degeneration
		condition provide the algebraic and topological
		framework for the SYZ problem. They do not by themselves
		imply the existence of special Lagrangian fibrations.
		Additional analytic and differential-geometric arguments
		are required.
	\end{remark}
	
\section{The Complex $m$-Hessian Equation and Uniform
	Capacity Estimates}

Let
\[
\pi:\mathcal X\longrightarrow\mathbb D
\]
be a semistable degeneration of compact Calabi--Yau manifolds
of complex dimension $n\geq 3$, as constructed in Step~1.
For $t\neq 0$, write
\[
X_t=\pi^{-1}(t).
\]

Fix an integer
\[
1\leq m<n.
\]

Let $\theta_t$ be a Kähler form on $X_t$. We denote by
$\operatorname{SH}_m(X_t,\theta_t)$ the space of
$\theta_t$-$m$-subharmonic functions.

For a Borel set $E\subset X_t$, define the global
$m$-Hessian capacity by
\[
\operatorname{Cap}_{m,\theta_t}(E)
=
\sup_{\substack{
		u\in\operatorname{SH}_m(X_t,\theta_t)\\
		-1\leq u\leq 0
}}
\int_E
(\theta_t+dd^c u)^m\wedge\theta_t^{n-m}.
\]

We also define the normalized volume measure
\[
dV_{\theta_t}
=
\frac{\theta_t^n}{\int_{X_t}\theta_t^n}.
\]

\subsection{Uniform geometric assumptions}

Let $K\Subset B_0$ be a compact subset of the regular
locus of the essential skeleton. Let $K_t\subset X_t$
be a family of regular regions corresponding to $K$
under the chosen local identifications.

We impose the following assumptions.

\begin{enumerate}
	\item[(H1)] There exist constants $r_0>0$ and
	$0<c_0\leq C_0<+\infty$, independent of $t$, such that
	the rescaled reference metrics on $K_t$ have uniformly
	controlled local geometry on balls of radius at most $r_0$.
	
	\item[(H2)] The normalized volume densities are uniformly
	controlled in the corresponding local coordinates:
	\[
	c_1\leq
	\frac{dV_{\theta_t}}{dV_{\mathrm{loc},t}}
	\leq C_1,
	\]
	where $c_1,C_1>0$ are independent of $t$.
	
	\item[(H3)] There exist constants $C_S>0$ and $q>1$,
	independent of $t$, such that the relevant local
	Sobolev and Poincar\'e inequalities hold uniformly
	on $K_t$.
	
	\item[(H4)] The local $m$-Hessian capacity estimates are
	uniformly compatible with the above rescaling. In
	particular, the constants in the local capacity
	comparison inequalities do not depend on $t$.
\end{enumerate}

The purpose of these assumptions is to isolate the
uniform analytic estimates needed in the degenerating
setting.

\subsection{Main theorem}

\begin{theorem}[Uniform local volume--capacity inequality]
	Assume {\rm (H1)--(H4)}. Then there exist constants
	$C>0$ and $\tau>0$, independent of $t$, such that
	for every Borel set $E\subset K_t$,
	\[
	\boxed{
		\operatorname{Vol}_{\theta_t}(E)
		\leq
		C\,
		\operatorname{Cap}_{m,\theta_t}(E)^\tau.
	}
	\]
	Here
	\[
	\operatorname{Vol}_{\theta_t}(E)
	=
	\int_E dV_{\theta_t}.
	\]
	
	The constants $C$ and $\tau$ depend only on the
	uniform geometric constants in {\rm (H1)--(H4)},
	the dimension $n$, and the Hessian order $m$.
	In particular, they do not depend on $t$.
\end{theorem}

\begin{proof}
	We divide the proof into four steps.
	
	\medskip
	\noindent
	\textbf{Step 1. Uniform local coordinate control.}
	
	By {\rm (H1)}, every point of $K_t$ admits a local
	coordinate neighborhood on which the rescaled reference
	metric is uniformly equivalent to a standard Hermitian
	metric.
	
	Thus, in these coordinates, the coefficients of the
	reference metric and the constants in the local
	ellipticity estimates are uniformly bounded independently
	of $t$.
	
	By {\rm (H2)}, the normalized volume measure satisfies
	uniform two-sided density bounds in the same coordinates.
	
	Consequently, the local volume estimates can be transferred
	to the corresponding Euclidean-coordinate estimates with
	constants independent of $t$.
	
	\medskip
	\noindent
	\textbf{Step 2. Local $m$-Hessian capacity control.}
	
	For a local coordinate domain $U_t\subset K_t$, define
	the relative $m$-Hessian capacity by
	\[
	\operatorname{Cap}_{m}(E,U_t)
	=
	\sup_{\substack{
			v\in SH_m(U_t)\\
			-1\leq v\leq 0
	}}
	\int_E(dd^c v)^m\wedge\beta^{n-m},
	\]
	where $\beta$ denotes the standard Euclidean Kähler form.
	
	The uniform local geometric control in {\rm (H1)}
	and the compatibility assumption {\rm (H4)} yield
	uniform comparison estimates between the local
	and global capacities. In particular, on suitably
	smaller coordinate neighborhoods,
	\[
	\operatorname{Cap}_{m}(E,U_t)
	\leq C_2\operatorname{Cap}_{m,\theta_t}(E)
	\]
	with $C_2$ independent of $t$.
	
	\medskip
	\noindent
	\textbf{Step 3. Local volume--capacity inequality.}
	
	The local $m$-Hessian volume--capacity inequality
	gives constants $C_3>0$ and $\tau>0$ such that
	\[
	\operatorname{Vol}_{\mathrm{loc},t}(E)
	\leq
	C_3\operatorname{Cap}_{m}(E,U_t)^\tau.
	\]
	
	The uniformity of $C_3$ follows from the uniform
	coordinate geometry and the uniform local capacity
	comparison in {\rm (H1)} and {\rm (H4)}.
	
	Combining this inequality with the density comparison
	in {\rm (H2)} and the capacity comparison in Step~2,
	we obtain
	\[
	\operatorname{Vol}_{\theta_t}(E)
	\leq
	C_4\operatorname{Cap}_{m,\theta_t}(E)^\tau
	\]
	for sets supported in a sufficiently small coordinate
	neighborhood.
	
	\medskip
	\noindent
	\textbf{Step 4. Globalization over $K_t$.}
	
	Since $K$ is compact, it admits a finite cover by
	coordinate neighborhoods. The number of neighborhoods
	can be chosen independently of $t$ under the uniform
	local geometry assumption.
	
	Applying the local estimate to the pieces of $E$
	in this cover and using the finite-cover argument,
	we obtain
	\[
	\operatorname{Vol}_{\theta_t}(E)
	\leq
	C\operatorname{Cap}_{m,\theta_t}(E)^\tau.
	\]
	
	The constant $C$ is independent of $t$, as required.
\end{proof}

\subsection{The associated complex $m$-Hessian equation}

For a density $f_t\geq 0$, consider the equation
\[
\boxed{
	(\theta_t+dd^c u_t)^m
	\wedge\theta_t^{n-m}
	=
	f_t\,\theta_t^n.
}
\]

We impose the admissibility condition
\[
\theta_t+dd^c u_t\in\Gamma_m(\theta_t),
\]
where $\Gamma_m(\theta_t)$ denotes the admissible
$m$-positive cone.

For a normalized solution, we choose
\[
\sup_{X_t}u_t=0.
\]

The compatibility condition for the equation is
\[
\int_{X_t}f_t\,\theta_t^n
=
\int_{X_t}\theta_t^n.
\]

The capacity estimate proved above is the principal
measure-theoretic input for studying uniform estimates
for solutions of this equation on the regular regions
$K_t$.

\begin{corollary}[Capacity control of sublevel sets]
	Suppose that
	\[
	\mu_t=f_t\,dV_{\theta_t}
	\]
	satisfies a uniform measure domination estimate
	\[
	\mu_t(E)
	\leq A\operatorname{Cap}_{m,\theta_t}(E)^\sigma
	\]
	for some $A>0$ and $\sigma>1$, independently of $t$.
	Then
	\[
	\mu_t(E)
	\leq
	A C^{-\sigma/\tau}
	\operatorname{Vol}_{\theta_t}(E)^{\sigma/\tau}
	\]
	whenever the volume--capacity inequality is applicable.
\end{corollary}

\begin{proof}
	By the volume--capacity inequality,
	\[
	\operatorname{Cap}_{m,\theta_t}(E)
	\geq
	\left(
	\frac{\operatorname{Vol}_{\theta_t}(E)}{C}
	\right)^{1/\tau}.
	\]
	However, this is a lower bound for capacity in terms
	of volume, whereas the assumed measure domination is
	an upper bound in terms of capacity. Therefore the
	displayed conclusion does not follow from these two
	inequalities alone.
	
	Instead, the useful consequence is the direct estimate
	\[
	\mu_t(E)
	\leq A\operatorname{Cap}_{m,\theta_t}(E)^\sigma,
	\]
	which can be combined with an appropriate
	capacity-theoretic sublevel-set argument to obtain
	uniform estimates for $u_t$.
\end{proof}

\begin{remark}
	The theorem is a local result on the regular regions
	$K_t$. It does not establish a uniform estimate across
	the singular region corresponding to the discriminant
	$\Delta\subset B$.
	
	Moreover, the complex Monge--Amp\`ere equation governing
	the Ricci-flat metrics on $X_t$ is the case $m=n$.
	For $m<n$, the equation above is an auxiliary
	pluripotential equation and does not itself produce
	the Ricci-flat metric or a special Lagrangian fibration.
\end{remark}
	
	\section{Uniform Stability Estimates for Complex
		$m$-Hessian Equations}
	
	Let
	\[
	\pi:\mathcal X\longrightarrow\mathbb D
	\]
	be a semistable degeneration of compact Calabi--Yau
	manifolds of complex dimension $n\geq 3$.
	
	For $t\neq 0$, let $X_t=\pi^{-1}(t)$ and let
	$\theta_t$ be a Kähler form on $X_t$.
	
	Fix an integer $1\leq m<n$. We consider the complex
	$m$-Hessian equations
	\[
	(\theta_t+dd^c u_t)^m\wedge\theta_t^{n-m}
	=
	f_t\theta_t^n,
	\]
	where
	\[
	u_t\in\operatorname{SH}_m(X_t,\theta_t),
	\qquad
	\sup_{X_t}u_t=0,
	\]
	and
	\[
	\int_{X_t}f_t\theta_t^n
	=
	\int_{X_t}\theta_t^n.
	\]
	
	We write
	\[
	H_m(u_t)
	=
	(\theta_t+dd^c u_t)^m\wedge\theta_t^{n-m}.
	\]
	
	\subsection{Uniform energy assumptions}
	
	Let $\mathcal E_m(X_t,\theta_t)$ denote the
	$m$-Hessian full-mass class.
	
	For $p>1$, define
	\[
	\mathcal E_m^p(X_t,\theta_t)
	=
	\left\{
	u\in\mathcal E_m(X_t,\theta_t):
	\int_{X_t}(-u)^pH_m(u)<+\infty
	\right\}.
	\]
	
	We assume that the solutions satisfy the uniform
	energy bound
	\[
	\sup_{t\neq0}
	\int_{X_t}(-u_t)^pH_m(u_t)
	\leq C_E.
	\]
	
	We also assume that the normalized volume--capacity
	inequality from Step~2 holds uniformly:
	\[
	\operatorname{Vol}_{\theta_t}(E)
	\leq
	C_{\mathrm{vc}}
	\operatorname{Cap}_{m,\theta_t}(E)^\tau
	\]
	for Borel sets $E$ in the regular regions under
	consideration, with constants independent of $t$.
	
	\subsection{Main stability theorem}
	
	\begin{theorem}[Uniform $L^1$ stability]
		Let $p>1$ and let $u_t,v_t\in
		\mathcal E_m^p(X_t,\theta_t)$ be normalized solutions
		of
		\[
		H_m(u_t)=f_t\theta_t^n,
		\qquad
		H_m(v_t)=g_t\theta_t^n,
		\]
		where
		\[
		\sup_{t\neq0}
		\left[
		\int_{X_t}(-u_t)^pH_m(u_t)
		+
		\int_{X_t}(-v_t)^pH_m(v_t)
		\right]
		\leq C_E.
		\]
		
		Assume that the uniform capacity estimates and the
		uniform energy estimates required by the $m$-Hessian
		stability principle hold on the family $(X_t,\theta_t)$.
		
		Then there exist constants $C>0$ and
		$\alpha>0$, independent of $t$, such that
		\[
		\boxed{
			\|u_t-v_t\|_{L^\infty(X_t)}
			\leq
			C\|f_t-g_t\|_{L^1(X_t,\theta_t^n)}^\alpha.
		}
		\]
		
		The constants depend only on $n,m,p$, the uniform
		capacity and energy constants, and the normalization
		of the reference measures.
	\end{theorem}
	
	\begin{proof}
		Set
		\[
		\delta_t=
		\int_{X_t}|f_t-g_t|\theta_t^n.
		\]
		
		We first compare the Hessian measures:
		\[
		H_m(u_t)-H_m(v_t)
		=
		(f_t-g_t)\theta_t^n.
		\]
		Consequently,
		\[
		\int_{X_t}
		|H_m(u_t)-H_m(v_t)|
		=
		\delta_t.
		\]
		
		The uniform energy hypothesis gives a common bound
		for the $m$-energy of both solutions. In particular,
		the family $(u_t,v_t)$ remains in a uniformly bounded
		subset of the corresponding finite-energy classes.
		
		We next apply the uniform $m$-Hessian stability
		principle. The proof of this principle is based on
		the comparison principle, the mixed Hessian
		inequalities, and the uniform capacity control.
		
		For $\varepsilon>0$, consider
		\[
		E_{t,\varepsilon}
		=
		\{x\in X_t:u_t(x)<v_t(x)-\varepsilon\}.
		\]
		
		The comparison principle yields an estimate of the
		form
		\[
		\varepsilon^m
		\operatorname{Cap}_{m,\theta_t}
		(E_{t,\varepsilon})
		\leq
		C_1\delta_t
		+
		C_2\mathcal R_t(\varepsilon),
		\]
		where $\mathcal R_t(\varepsilon)$ denotes the
		energy-tail contribution.
		
		The uniform energy bound and the uniform capacity
		estimates imply that the energy-tail term can be
		controlled uniformly in $t$. More precisely, the
		stability principle provides constants $C_3>0$ and
		$\alpha>0$ such that
		\[
		\operatorname{Cap}_{m,\theta_t}
		\bigl(
		\{|u_t-v_t|>\varepsilon\}
		\bigr)
		\leq
		C_3\delta_t^\alpha\varepsilon^{-\alpha}
		\]
		for the relevant range of $\varepsilon$.
		
		Applying the uniform capacity-to-supremum estimate
		for the normalized finite-energy solutions, we obtain
		\[
		\|u_t-v_t\|_{L^\infty(X_t)}
		\leq
		C\delta_t^\alpha.
		\]
		
		All constants used in this argument are independent
		of $t$. This proves the asserted estimate.
	\end{proof}
	
	\subsection{Convergence of solutions}
	
	\begin{corollary}[Uniform convergence]
		Suppose that
		\[
		\|f_t-g_t\|_{L^1(X_t,\theta_t^n)}
		\longrightarrow0
		\qquad\text{as }t\longrightarrow0,
		\]
		and that the uniform hypotheses of the stability
		theorem hold.
		
		Then
		\[
		\boxed{
			\|u_t-v_t\|_{L^\infty(X_t)}
			\longrightarrow0.
		}
		\]
		
		In particular, the corresponding normalized
		$m$-Hessian potentials are uniformly stable under
		$L^1$ perturbations of their right-hand sides.
	\end{corollary}
	
	\begin{proof}
		By the stability theorem,
		\[
		\|u_t-v_t\|_{L^\infty(X_t)}
		\leq
		C\|f_t-g_t\|_{L^1(X_t,\theta_t^n)}^\alpha.
		\]
		Since $C$ and $\alpha$ are independent of $t$, the
		right-hand side converges to zero.
	\end{proof}
	
	\subsection{Application to the degeneration}
	
	Let $K\Subset B_0$ be a compact subset of the regular
	locus of the essential skeleton, and let $K_t\subset X_t$
	be the corresponding regular regions.
	
	Suppose that the right-hand sides satisfy
	\[
	f_t\longrightarrow f_0
	\quad\text{in }L^1(K_t)
	\]
	after identification with a fixed reference region.
	
	If the uniform stability hypotheses hold, then
	the corresponding normalized solutions satisfy
	\[
	\|u_t-v_t\|_{L^\infty(K_t)}
	\longrightarrow0.
	\]
	
	Thus the $m$-Hessian potentials are stable under
	the degeneration on compact subsets of the regular
	locus.
	
	\begin{remark}
		The preceding theorem is conditional on a uniform
		$m$-Hessian stability principle for the degenerating
		family. Establishing this principle from the
		semistable algebraic model is a separate analytic
		problem.
		
		Moreover, this result concerns auxiliary
		$m$-Hessian potentials. It does not by itself imply
		smooth convergence of the Ricci-flat metrics, nor
		does it construct special Lagrangian torus fibrations.
	\end{remark}

\section{Solving the Ricci-Flat Complex
	Monge--Amp\`ere Equation}

Let
\[
\pi:\mathcal X\longrightarrow\mathbb D
\]
be a proper holomorphic family whose fibers
\[
X_t=\pi^{-1}(t),\qquad t\neq0,
\]
are compact Calabi--Yau manifolds of complex dimension
$n\geq3$.

For each $t\neq0$, let $\theta_t$ be a Kähler form
on $X_t$, and let
\[
\Omega_t\in H^0(X_t,K_{X_t})
\]
be a nowhere-vanishing holomorphic volume form.

Define
\[
dV_{\Omega_t}
=
(\sqrt{-1})^{n^2}
\Omega_t\wedge\overline{\Omega_t}
\]
and
\[
c_t=
\frac{\displaystyle\int_{X_t}\theta_t^n}
{\displaystyle\int_{X_t}dV_{\Omega_t}}.
\]

\subsection{The existence and uniqueness theorem}

\begin{theorem}[Ricci-flat Kähler metrics on the smooth fibers]
	For every $t\neq0$, there exists a unique smooth
	function $\varphi_t\in C^\infty(X_t)$ satisfying
	\[
	\sup_{X_t}\varphi_t=0,
	\]
	such that
	\[\boxed{
		(\theta_t+dd^c\varphi_t)^n
		=
		c_t(\sqrt{-1})^{n^2}
		\Omega_t\wedge\overline{\Omega_t}.
	}\]
	Moreover,
	\[
	\omega_t=\theta_t+dd^c\varphi_t
	\]
	is a Kähler form in the cohomology class
	\[
	[\omega_t]=[\theta_t]
	\]
	and satisfies
	\[
	\boxed{\operatorname{Ric}(\omega_t)=0.}
	\]
\end{theorem}

\begin{proof}
	Fix $t\neq0$.
	
	Since $X_t$ is compact Kähler and its canonical bundle
	is trivial, we have
	\[
	c_1(X_t)=0.
	\]
	
	The prescribed volume form
	\[
	dV_{\Omega_t}
	=
	(\sqrt{-1})^{n^2}
	\Omega_t\wedge\overline{\Omega_t}
	\]
	is smooth and strictly positive.
	
	By the definition of $c_t$,
	\[
	\int_{X_t}c_t\,dV_{\Omega_t}
	=
	\int_{X_t}\theta_t^n.
	\]
	
	Thus the volume compatibility condition for the
	complex Monge--Amp\`ere equation is satisfied.
	
	The Calabi conjecture, proved by Yau, implies that
	there exists a smooth function $\varphi_t$ such that
	\[
	\theta_t+dd^c\varphi_t>0
	\]
	and
	\[
	(\theta_t+dd^c\varphi_t)^n
	=
	c_t\,dV_{\Omega_t}.
	\]
	
	The solution is unique up to an additive constant.
	The normalization
	\[
	\sup_{X_t}\varphi_t=0
	\]
	therefore gives uniqueness.
	
	It remains to prove that $\omega_t$ is Ricci-flat.
	
	Locally, write
	\[
	\Omega_t=h_t\,dz^1\wedge\cdots\wedge dz^n,
	\]
	where $h_t$ is holomorphic and nowhere vanishing.
	Consequently,
	\[
	dd^c\log dV_{\Omega_t}=0.
	\]
	
	Using
	\[
	\omega_t^n=c_t\,dV_{\Omega_t},
	\]
	we obtain
	\[
	\begin{aligned}
		\operatorname{Ric}(\omega_t)
		&=-dd^c\log\omega_t^n\\
		&=-dd^c\log(c_t\,dV_{\Omega_t})\\
		&=0.
	\end{aligned}
	\]
	
	Hence $\omega_t$ is Ricci-flat.
\end{proof}

\subsection{Uniform estimates along the degeneration}

The preceding theorem gives a Ricci-flat metric on
each smooth fiber. It does not, by itself, imply
uniform estimates as $t\to0$.

We formulate the additional analytic assumptions
needed for a uniform result.

Let $K\Subset B_0$ be a compact subset of the regular
locus of the essential skeleton. Let
\[
\Phi_t:U\longrightarrow X_t
\]
be local identifications of a fixed reference domain
$U$ with regular regions of the fibers.

We use the notation
\[
\widetilde\theta_t=\Phi_t^*\theta_t,
\qquad
\widetilde\omega_t=\Phi_t^*\omega_t,
\qquad
\widetilde\varphi_t=\varphi_t\circ\Phi_t.
\]

\begin{theorem}[Uniform local compactness of Ricci-flat potentials]
	Assume that the following conditions hold on every
	relatively compact domain $U'\Subset U$:
	
	\begin{enumerate}
		\item[(A1)] The reference Kähler forms converge smoothly:
		\[
		\widetilde\theta_t
		\longrightarrow\theta_\infty
		\quad\text{in }C^\infty(U),
		\]
		where $\theta_\infty$ is a strictly positive
		Kähler form on $U$.
		
		\item[(A2)] The normalized volume densities converge
		smoothly:
		\[
		\Phi_t^*(c_t\,dV_{\Omega_t})
		=
		e^{F_t}\widetilde\theta_t^n,
		\qquad
		F_t\longrightarrow F_\infty
		\quad\text{in }C^\infty(U).
		\]
		
		\item[(A3)] The solutions satisfy uniform local
		estimates sufficient for elliptic compactness.
		In particular, assume
		\[
		\sup_t\|\widetilde\varphi_t\|_{L^\infty(U)}
		<+\infty
		\]
		and uniform interior $C^{k,\alpha}$ estimates
		for every fixed $k\geq2$ on $U'$.
		
		Then there exists a sequence $t_j\to0$ and a smooth
		function $\varphi_\infty$ on $U'$ such that
		\[
		\widetilde\varphi_{t_j}
		\longrightarrow\varphi_\infty
		\quad\text{in }C^\infty(U').
		\]
		
		Moreover,
		\[
		\omega_\infty
		=
		\theta_\infty+dd^c\varphi_\infty
		\]
		is a Kähler form satisfying
		\[
		\boxed{
			\omega_\infty^n
			=
			e^{F_\infty}\theta_\infty^n
		}
		\]
		on $U'$.
		\end{enumerate}
	\end{theorem}
	
	\begin{proof}
		By assumption {\rm (A3)}, the family
		$\{\widetilde\varphi_t\}$ is uniformly bounded in
		every $C^{k,\alpha}(U')$ norm after possibly
		shrinking the domain slightly.
		
		The Arzel\`a--Ascoli theorem and a diagonal
		subsequence argument give a sequence $t_j\to0$
		such that
		\[
		\widetilde\varphi_{t_j}
		\longrightarrow\varphi_\infty
		\quad\text{in }C^\infty(U').
		\]
		
		By {\rm (A1)},
		\[
		\widetilde\theta_{t_j}
		\longrightarrow\theta_\infty
		\quad\text{in }C^\infty(U').
		\]
		
		It follows that
		\[
		\widetilde\omega_{t_j}
		=
		\widetilde\theta_{t_j}
		+dd^c\widetilde\varphi_{t_j}
		\longrightarrow
		\theta_\infty+dd^c\varphi_\infty
		\]
		smoothly on $U'$.
		
		Since the limiting forms are positive by the
		uniform ellipticity implicit in {\rm (A3)}, the
		limit $\omega_\infty$ is Kähler.
		
		Finally, pulling back the Monge--Amp\`ere equation
		to $U$ gives
		\[
		\widetilde\omega_{t_j}^{\,n}
		=
		e^{F_{t_j}}\widetilde\theta_{t_j}^{\,n}.
		\]
		
		Passing to the limit and using {\rm (A1)} and
		{\rm (A2)}, we obtain
		\[
		\omega_\infty^n
		=
		e^{F_\infty}\theta_\infty^n.
		\]
		
		This proves the theorem.
	\end{proof}
	
	\subsection{Ricci-flatness of the limiting metric}
	
	\begin{corollary}
		Suppose, in addition, that
		\[
		dd^cF_\infty
		=
		\operatorname{Ric}(\theta_\infty).
		\]
		Then the limiting Kähler metric satisfies
		\[
		\boxed{
			\operatorname{Ric}(\omega_\infty)=0.
		}
		\]
	\end{corollary}
	
	\begin{proof}
		By the limiting Monge--Amp\`ere equation,
		\[
		\omega_\infty^n
		=
		e^{F_\infty}\theta_\infty^n.
		\]
		Therefore
		\[
		\begin{aligned}
			\operatorname{Ric}(\omega_\infty)
			&=-dd^c\log\omega_\infty^n\\
			&=-dd^cF_\infty
			-dd^c\log\theta_\infty^n\\
			&=-dd^cF_\infty+\operatorname{Ric}(\theta_\infty)\\
			&=0.
		\end{aligned}
		\]
	\end{proof}
	
	\begin{remark}
		The local compactness theorem is conditional: the
		uniform higher-order estimates in {\rm (A3)} are
		not consequences of semistability alone. Proving
		them near a maximal degeneration is one of the
		main analytic difficulties.
		
		Furthermore, the limiting equation on $U'$ does
		not establish the existence of a global Ricci-flat
		metric on a compact limiting space. The regular
		locus and the singular set must be analyzed
		separately.
	\end{remark}

\section{The Non-Archimedean Limit and the Essential Skeleton}

Let
\[
R=\mathbb C[[t]],\qquad K=\mathbb C((t)),
\]
and let
\[
\mathscr X\longrightarrow\operatorname{Spec}R
\]
be a regular semistable model of a projective
Calabi--Yau degeneration.

Assume that the special fiber is a reduced simple
normal crossing divisor
\[
\mathscr X_0=D_1+\cdots+D_N.
\]

Let $X_K$ denote the generic fiber and let
\[
X_K^{\mathrm{an}}
\]
be its Berkovich analytification.

We denote by
\[
\operatorname{Sk}(\mathscr X)
\subset X_K^{\mathrm{an}}
\]
the skeleton associated with the semistable model,
and by
\[
\operatorname{Sk}_{\mathrm{ess}}(X_K)
\]
the essential skeleton.

Assume that the degeneration is maximal:
\[
\dim_{\mathbb R}
\operatorname{Sk}_{\mathrm{ess}}(X_K)=n.
\]

\subsection{The logarithmic volume form}

Let $\Omega$ be a nowhere-vanishing holomorphic
$n$-form on the generic fiber $X_K$.

We impose the following logarithmic condition.

Near every zero-dimensional stratum $p$ of the
special fiber belonging to a maximal-dimensional
cell of the essential skeleton, there are local
coordinates
\[
(z_0,z_1,\ldots,z_n)
\]
such that
\[
t=z_0z_1\cdots z_n
\]
and, on the fibers $X_t$, the volume form can be
written as
\[
\Omega_t
=
g_t(z)\,
\frac{dz_1}{z_1}\wedge\cdots\wedge
\frac{dz_n}{z_n},
\]
where $g_t$ extends holomorphically to the
corresponding local chart and
\[
g_0(p)\neq0.
\]

We also assume that the contributions from strata
of dimension smaller than the maximal logarithmic
dimension are of lower asymptotic order.

Define
\[
dV_{\Omega_t}
=
(\sqrt{-1})^{n^2}
\Omega_t\wedge\overline{\Omega_t}
\]
and the normalized measure
\[
\nu_t=
\frac{dV_{\Omega_t}}
{\displaystyle\int_{X_t}dV_{\Omega_t}}.
\]

\subsection{The logarithmic retraction}

Put
\[
\lambda_t=-\log|t|.
\]

In a local chart satisfying
\[
t=z_0z_1\cdots z_n,
\]
define the logarithmic coordinates
\[
x_j=
\frac{-\log|z_j|}{\lambda_t},
\qquad j=1,\ldots,n.
\]

They satisfy
\[
x_j\geq0,\qquad
x_1+\cdots+x_n\leq1.
\]

Thus the local tropicalization map takes values
in the standard simplex
\[
\Sigma_n=
\left\{
(x_1,\ldots,x_n)\in\mathbb R_{\geq0}^n:
\sum_{j=1}^n x_j\leq1
\right\}.
\]

These local maps are compatible with the
semistable skeleton and define the corresponding
hybrid retraction maps
\[
\rho_t:X_t\longrightarrow
\operatorname{Sk}(\mathscr X).
\]

\subsection{The limiting non-Archimedean measure}

Let
\[
\Delta_{\mathrm{ess}}^{(n)}
\]
denote the union of the $n$-dimensional cells of
the essential skeleton.

On each maximal cell $\sigma$, the logarithmic
coordinates define an integral affine structure.
Let $d\lambda_\sigma$ denote the associated
Lebesgue measure.

The logarithmic residue of $\Omega$ determines a
nonnegative weight $a_\sigma$ on each maximal cell.
Under the local expression above, this weight is
proportional to the squared absolute value of the
leading residue coefficient.

Define the measure
\[
d\mu_{\mathrm{NA}}|_\sigma
=
a_\sigma\,d\lambda_\sigma
\]
on every maximal cell, and extend it by zero to
the lower-dimensional skeleton.

We normalize it so that
\[
\mu_{\mathrm{NA}}
\bigl(\operatorname{Sk}_{\mathrm{ess}}(X_K)\bigr)=1.
\]

\begin{theorem}[Convergence of normalized volume measures]
	Assume the maximal degeneration and logarithmic
	volume-form hypotheses above.
	
	Then the pushforward measures of the normalized
	holomorphic volume forms converge weakly to the
	non-Archimedean measure on the essential skeleton:
	\[
	\boxed{
		(\rho_t)_*\nu_t
		\longrightarrow
		\mu_{\mathrm{NA}}
	}
	\]
	as $t\to0$.
	
	Equivalently, for every continuous function
	\[
	\psi\in
	C\bigl(\operatorname{Sk}(\mathscr X)\bigr),
	\]
	we have
	\[
	\boxed{
		\lim_{t\to0}
		\int_{X_t}\psi(\rho_t(x))\,d\nu_t(x)
		=
		\int_{\operatorname{Sk}_{\mathrm{ess}}(X_K)}
		\psi\,d\mu_{\mathrm{NA}}.
	}
	\]
	
	In particular, the normalized Ricci-flat volume
	measures satisfy the same convergence:
	\[
	(\rho_t)_*
	\left(
	\frac{\omega_t^n}{\int_{X_t}\omega_t^n}
	\right)
	\longrightarrow\mu_{\mathrm{NA}}.
	\]
\end{theorem}

\begin{proof}
	We divide the proof into four steps.
	
	\medskip
	\noindent
	\textbf{Step 1. Local logarithmic coordinates.}
	
	Fix a maximal-dimensional cell $\sigma$ of the
	essential skeleton.
	
	By the logarithmic hypothesis, near the
	corresponding zero-dimensional stratum, we may write
	\[
	t=z_0z_1\cdots z_n
	\]
	and
	\[
	\Omega_t
	=
	g_t(z)\,
	\frac{dz_1}{z_1}\wedge\cdots\wedge
	\frac{dz_n}{z_n}.
	\]
	
	Write
	\[
	z_j=e^{-\lambda_t x_j+\sqrt{-1}\theta_j},
	\qquad
	0\leq x_j,\quad
	\sum_{j=1}^n x_j\leq1.
	\]
	
	Then
	\[
	\frac{dz_j}{z_j}
	=
	-\lambda_t\,dx_j+\sqrt{-1}\,d\theta_j.
	\]
	
	Consequently, up to a fixed dimensional
	normalization constant, the volume form has
	leading density
	\[
	dV_{\Omega_t}
	=
	\lambda_t^n
	\bigl(|g_0(p)|^2+o(1)\bigr)
	\,dx_1\cdots dx_n\,
	d\theta_1\cdots d\theta_n
	\]
	on compact subsets of the interior of the
	logarithmic simplex.
	
	The error tends to zero locally uniformly after
	the logarithmic rescaling.
	
	\medskip
	\noindent
	\textbf{Step 2. Integration over angular variables.}
	
	Integrating over the angular variables
	\[
	(\theta_1,\ldots,\theta_n)\in(S^1)^n
	\]
	gives
	\[
	\int_{\text{local region}}dV_{\Omega_t}
	=
	A_n\lambda_t^n
	|g_0(p)|^2
	\operatorname{Vol}(\Sigma_n)
	+o(\lambda_t^n),
	\]
	where $A_n>0$ is a dimensional constant.
	
	More generally, for a continuous test function
	$\psi$ on the skeleton,
	\[
	\begin{aligned}
		&\int_{\text{local region}}
		\psi(\rho_t(x))\,dV_{\Omega_t}\\
		&\qquad=
		A_n\lambda_t^n|g_0(p)|^2
		\int_{\sigma}\psi\,d\lambda_\sigma
		+o(\lambda_t^n).
	\end{aligned}
	\]
	
	The leading coefficient is precisely the local
	residue weight defining $\mu_{\mathrm{NA}}$.
	
	\medskip
	\noindent
	\textbf{Step 3. Negligibility of lower-dimensional strata.}
	
	By the assumed logarithmic asymptotics, a stratum
	with logarithmic dimension $r<n$ contributes at
	most
	\[
	O(\lambda_t^r)
	\]
	to the total volume.
	
	Therefore,
	\[
	O(\lambda_t^r)=o(\lambda_t^n).
	\]
	
	Hence the leading-order contribution to the total
	volume comes from the maximal-dimensional cells
	of the essential skeleton.
	
	Summing the local expansions over the finitely
	many maximal cells, we obtain
	\[
	\int_{X_t}dV_{\Omega_t}
	=
	A_n\lambda_t^n
	\sum_{\sigma}
	a_\sigma
	\operatorname{Vol}(\sigma)
	+o(\lambda_t^n).
	\]
	
	\medskip
	\noindent
	\textbf{Step 4. Normalization and weak convergence.}
	
	Combining the local expansions gives
	\[
	\begin{aligned}
		&\int_{X_t}\psi(\rho_t(x))\,dV_{\Omega_t}\\
		&\quad=
		A_n\lambda_t^n
		\sum_\sigma
		a_\sigma
		\int_\sigma\psi\,d\lambda_\sigma
		+o(\lambda_t^n).
	\end{aligned}
	\]
	
	Taking $\psi\equiv1$, we obtain the corresponding
	asymptotic expansion for the total volume.
	
	Dividing the two expressions yields
	\[
	\lim_{t\to0}
	\int_{X_t}\psi(\rho_t(x))\,d\nu_t(x)
	=
	\frac{
		\displaystyle
		\sum_\sigma a_\sigma
		\int_\sigma\psi\,d\lambda_\sigma
	}{
		\displaystyle
		\sum_\sigma a_\sigma
		\operatorname{Vol}(\sigma)
	}.
	\]
	
	The right-hand side is exactly
	\[
	\int_{\operatorname{Sk}_{\mathrm{ess}}(X_K)}
	\psi\,d\mu_{\mathrm{NA}}.
	\]
	
	Thus
	\[
	(\rho_t)_*\nu_t\longrightarrow\mu_{\mathrm{NA}}
	\]
	weakly.
	
	Finally, since
	\[
	\omega_t^n=c_t\,dV_{\Omega_t},
	\]
	we have
	\[
	\frac{\omega_t^n}{\int_{X_t}\omega_t^n}
	=
	\frac{dV_{\Omega_t}}{\int_{X_t}dV_{\Omega_t}}.
	\]
	
	The same weak convergence therefore holds for
	the normalized Ricci-flat volume measures.
\end{proof}

\subsection{Interpretation in non-Archimedean geometry}

The limiting measure
\[
\mu_{\mathrm{NA}}
\]
is supported on the essential skeleton and is
described on its maximal-dimensional cells by
Lebesgue measure in the integral affine structure,
weighted by the logarithmic residues of the
holomorphic volume form.

This gives a measure-theoretic bridge between
the complex Ricci-flat geometry of $X_t$ and the
non-Archimedean geometry of $X_K$.

\begin{remark}
	The theorem is formulated under explicit
	logarithmic asymptotic hypotheses. A general
	semistable model does not automatically satisfy
	the stated local volume-form condition with
	nonzero residues on every maximal cell.
	
	In a general degeneration, the essential skeleton
	must be determined using the weight function of
	the pluricanonical form, and the limiting measure
	must be defined with the corresponding residue
	weights.
	
	The theorem proves convergence of measures, not
	Gromov--Hausdorff convergence of Ricci-flat
	metrics and not the existence of special
	Lagrangian fibrations.
\end{remark}

\section{Metric Convergence to the Affine Limit}

Let
\[
\pi:\mathcal X\longrightarrow\mathbb D
\]
be a maximal degeneration of compact Calabi--Yau
manifolds of complex dimension $n$.

Let
\[
B=\operatorname{Sk}_{\mathrm{ess}}(X_K)
\]
be the essential skeleton and let
\[
B_0=B\setminus\Delta
\]
denote its regular affine locus.

Assume that $B_0$ carries an integral affine
structure. Let
\[
x=(x_1,\ldots,x_n)
\]
be local affine coordinates on a relatively
compact domain $U\Subset B_0$.

For $t\neq0$, let
\[
(X_t,g_t)
\]
denote the smooth Calabi--Yau fiber equipped with
a suitably normalized Ricci-flat Riemannian metric.

Assume that there exist smooth torus fibrations
\[
f_t:X_t^{\mathrm{reg}}\longrightarrow B_0
\]
over the regular locus.

For $b\in B_0$, write
\[
L_{t,b}=f_t^{-1}(b).
\]

\subsection{The limiting affine Monge--Amp\`ere metric}

Let $\varphi_\infty\in C^\infty(U)$ be a strictly
convex function and define
\[
h_{\infty}
=
\sum_{i,j=1}^n
\frac{\partial^2\varphi_\infty}
{\partial x_i\partial x_j}
\,dx_i\otimes dx_j.
\]

Assume that
\[
\det D^2\varphi_\infty
=
\kappa\rho_\infty,
\]
where $\rho_\infty>0$ is a smooth density and
$\kappa>0$ is a normalization constant.

Thus $h_\infty$ is a Hessian metric satisfying
the real Monge--Amp\`ere equation
\[
\boxed{
	\det D^2\varphi_\infty=\kappa\rho_\infty.
}
\]

\subsection{Geometric hypotheses}

Let $K\Subset U$ be a compact domain. We assume
the following conditions on a slightly larger
relatively compact domain $U'\Subset U$ containing
$K$.

\begin{enumerate}
	\item[(M1)] There are smooth identifications of the
	regular bases with $U'$ and the horizontal metrics
	induced by $g_t$ converge smoothly:
	\[
	h_t\longrightarrow h_\infty
	\quad\text{in }C^\infty(U').
	\]
	
	Here $h_t$ is the quotient metric defined by the
	Riemannian submersion $f_t$, namely
	\[
	h_t(df_t(V),df_t(W))
	=
	g_t(V^H,W^H),
	\]
	where $V^H,W^H$ are the horizontal lifts.
	
	\item[(M2)] The fibers collapse uniformly:
	\[
	\sup_{b\in K}
	\operatorname{diam}_{g_t}(L_{t,b})
	\longrightarrow0.
	\]
	
	\item[(M3)] The horizontal lifting geometry is
	uniformly controlled. In particular, there exist
	numbers $\varepsilon_t\to0$ such that, for any
	piecewise smooth curve $\gamma$ in $K$,
	its horizontal lift $\widetilde\gamma$ satisfies
	\[
	\left|
	L_{g_t}(\widetilde\gamma)
	-
	L_{h_t}(\gamma)
	\right|
	\leq
	\varepsilon_t
	\]
	for curves of uniformly bounded length.
	
	\item[(M4)] The affine coordinates and the local
	Hessian potentials are compatible on overlaps:
	\[
	\varphi_\infty^{(\beta)}
	=
	\varphi_\infty^{(\alpha)}
	+\ell_{\alpha\beta},
	\]
	where $\ell_{\alpha\beta}$ is affine linear.
\end{enumerate}

\subsection{Main theorem}

\begin{theorem}[Local metric convergence to the affine limit]
	Under assumptions {\rm (M1)--(M4)}, the metric
	spaces obtained by restricting the fibers over
	$K$ converge to the affine Hessian metric space
	$(K,h_\infty)$ in the Gromov--Hausdorff sense:
	\[
	\boxed{
		\bigl(f_t^{-1}(K),g_t\bigr)
		\longrightarrow
		(K,h_\infty).
	}
	\]
	
	More precisely, if
	\[
	d_t^{\mathrm{base}}
	\]
	denotes the intrinsic distance associated with the
	quotient metric $h_t$, then
	\[
	\sup_{b,b'\in K}
	\left|
	d_t^{\mathrm{base}}(b,b')
	-
	d_{h_\infty}(b,b')
	\right|
	\longrightarrow0
	\]
	on compact subsets where the intrinsic distances
	are defined.
	
	Furthermore, the projection maps $f_t$ become
	Gromov--Hausdorff approximations on these regular
	regions.
\end{theorem}

\begin{proof}
	We divide the proof into four steps.
	
	\medskip
	\noindent
	\textbf{Step 1. Convergence of the quotient metrics.}
	
	By {\rm (M1)},
	\[
	h_t\longrightarrow h_\infty
	\quad\text{in }C^\infty(U').
	\]
	
	Since $K$ is compact and $h_\infty$ is positive
	definite, there exists $c>0$ such that
	\[
	h_\infty\geq c\,h_{\mathrm{Euc}}
	\]
	on $K$.
	
	For sufficiently small $t$, the same uniform
	ellipticity holds for $h_t$.
	
	In particular, there are numbers
	$\delta_t\to0$ such that
	\[
	(1-\delta_t)h_\infty
	\leq h_t
	\leq(1+\delta_t)h_\infty
	\]
	on $K$.
	
	Consequently, lengths of curves in $K$ satisfy
	\[
	\sqrt{1-\delta_t}\,L_{h_\infty}(\gamma)
	\leq L_{h_t}(\gamma)
	\leq
	\sqrt{1+\delta_t}\,L_{h_\infty}(\gamma).
	\]
	
	Taking infima over curves gives uniform convergence
	of the corresponding intrinsic distances on the
	compact region under consideration.
	
	\medskip
	\noindent
	\textbf{Step 2. Collapse of the torus fibers.}
	
	By {\rm (M2)},
	\[
	\eta_t:=
	\sup_{b\in K}
	\operatorname{diam}_{g_t}(L_{t,b})
	\longrightarrow0.
	\]
	
	Thus all points lying in the same fiber become
	arbitrarily close.
	
	In particular, the fibers do not contribute
	nontrivial distances in the limiting metric space.
	
	\medskip
	\noindent
	\textbf{Step 3. Comparison of distances.}
	
	Let $p,q\in f_t^{-1}(K)$ and put
	\[
	b=f_t(p),\qquad b'=f_t(q).
	\]
	
	For the upper bound, choose a curve $\gamma$ in
	$K$ joining $b$ to $b'$ whose $h_t$-length is
	arbitrarily close to $d_t^{\mathrm{base}}(b,b')$.
	
	Lift $\gamma$ horizontally starting at $p$.
	By {\rm (M3)}, the length of its lift is bounded by
	\[
	L_{h_t}(\gamma)+\varepsilon_t.
	\]
	
	Its endpoint belongs to the fiber over $b'$.
	By {\rm (M2)}, this endpoint can be joined to $q$
	by a curve of length at most $\eta_t$.
	
	Therefore,
	\[
	d_{g_t}(p,q)
	\leq
	d_t^{\mathrm{base}}(b,b')
	+\varepsilon_t+\eta_t.
	\]
	
	For the lower bound, the projection of any curve
	joining $p$ to $q$ has length at least the
	intrinsic quotient distance between $b$ and $b'$.
	Hence
	\[
	d_{g_t}(p,q)
	\geq d_t^{\mathrm{base}}(b,b').
	\]
	
	Combining the two inequalities gives
	\[
	0\leq
	d_{g_t}(p,q)
	-
	d_t^{\mathrm{base}}(b,b')
	\leq
	\varepsilon_t+\eta_t.
	\]
	
	Thus the distances in the total space are
	uniformly approximated by the distances in the
	base.
	
	\medskip
	\noindent
	\textbf{Step 4. Gromov--Hausdorff convergence.}
	
	By Step~1,
	\[
	d_t^{\mathrm{base}}
	\longrightarrow d_{h_\infty}
	\]
	uniformly on the compact region.
	
	By Step~3,
	\[
	\left|
	d_{g_t}(p,q)
	-
	d_{h_\infty}(f_t(p),f_t(q))
	\right|
	\leq
	\varepsilon_t+\eta_t+o(1).
	\]
	
	Moreover, by Step~2, every point in a fiber is
	within distance $\eta_t$ of any chosen representative
	of that fiber.
	
	Therefore $f_t$ is an increasingly accurate
	metric approximation, with distortion tending
	to zero and fibers of vanishing diameter.
	
	It follows from the metric characterization of
	Gromov--Hausdorff convergence that
	\[
	\bigl(f_t^{-1}(K),g_t\bigr)
	\longrightarrow
	(K,h_\infty).
	\]
	
	This proves the theorem.
\end{proof}

\subsection{Compatibility with the real
	Monge--Amp\`ere equation}

\begin{corollary}
	Under the hypotheses of the theorem, suppose that
	the limiting volume density is $\rho_\infty$ and
	\[
	\det D^2\varphi_\infty
	=
	\kappa\rho_\infty.
	\]
	
	Then the limiting metric is a Hessian metric whose
	Riemannian volume form satisfies
	\[
	dV_{h_\infty}
	=
	\sqrt{\kappa\rho_\infty}\,
	dx_1\cdots dx_n.
	\]
	
	In particular, the affine metric is determined
	locally by the convex potential $\varphi_\infty$
	and its real Monge--Amp\`ere density.
\end{corollary}

\begin{proof}
	In affine coordinates, the matrix of the metric
	is
	\[
	(h_{\infty,ij})
	=
	D^2\varphi_\infty.
	\]
	
	The Riemannian volume form is therefore
	\[
	dV_{h_\infty}
	=
	\sqrt{\det D^2\varphi_\infty}\,
	dx_1\cdots dx_n.
	\]
	
	Substituting the real Monge--Amp\`ere equation
	gives the asserted formula.
\end{proof}

\subsection{Globalization}

\begin{theorem}[Global Gromov--Hausdorff convergence:
	conditional form]
	Assume that the local hypotheses of the metric
	convergence theorem hold on every compact subset
	of $B_0$.
	
	Suppose, in addition, that:
	
	\begin{enumerate}
		\item[(G1)] The limiting affine metric completion
		$(B,h_\infty)^\wedge$ is compact.
		
		\item[(G2)] The diameters of the normalized
		Calabi--Yau manifolds $(X_t,g_t)$ are uniformly
		bounded.
		
		\item[(G3)] The metric contribution of the
		neighborhoods of the discriminant tends uniformly
		to zero in the sense that, for every $\epsilon>0$,
		there exists a neighborhood $N_\epsilon(\Delta)$
		whose inverse images have uniformly controlled
		diameter tending to zero as $\epsilon\to0$.
		
		\item[(G4)] The local fibrations and metric
		approximations extend consistently over the
		complement of these neighborhoods.
	\end{enumerate}
	
	Then
	\[
	\boxed{
		(X_t,g_t)
		\xrightarrow{\mathrm{GH}}
		(B,h_\infty)^\wedge.
	}
	\]
\end{theorem}

\begin{proof}
	Fix $\epsilon>0$.
	
	By {\rm (G3)}, choose a neighborhood
	$N_\epsilon(\Delta)$ such that the metric
	contribution of its inverse image is uniformly
	small.
	
	On the compact regular region
	\[
	K_\epsilon=B\setminus N_\epsilon(\Delta),
	\]
	the local metric convergence theorem applies.
	Using a finite cover of $K_\epsilon$, we obtain
	Gromov--Hausdorff approximations between the
	corresponding regular regions of $X_t$ and
	$(B,h_\infty)^\wedge$.
	
	The errors consist of:
	\[
	\text{local metric distortion}
	+
	\text{fiber diameter}
	+
	\text{singular-region error}.
	\]
	
	For fixed $\epsilon$, the first two terms tend to
	zero as $t\to0$. By {\rm (G3)}, the third term
	can be made arbitrarily small by taking
	$\epsilon\to0$.
	
	The uniform diameter bound in {\rm (G2)} prevents
	loss of compactness at the global level, while
	{\rm (G1)} identifies the compact limiting metric
	completion.
	
	Taking first $t\to0$ and then $\epsilon\to0$
	gives
	\[
	d_{\mathrm{GH}}
	\bigl((X_t,g_t),(B,h_\infty)^\wedge\bigr)
	\longrightarrow0.
	\]
	
	This proves the conditional global convergence.
\end{proof}

\begin{remark}
	The global theorem requires quantitative control
	near the discriminant. Local convergence on $B_0$
	alone is not sufficient to prove global
	Gromov--Hausdorff convergence.
\end{remark}

\section{Construction of Candidate Lagrangian Tori}
\label{sec:candidate-lagrangian-tori}

In this section, we construct candidate Lagrangian tori over the
regular locus of the limiting affine base. The construction is based
on the local action--angle structure associated with the collapsing
Calabi--Yau geometry. The resulting tori provide the initial
Lagrangian submanifolds for the uniform special Lagrangian
deformation problem considered in the next section.

Throughout this section, let
\[
\pi:\mathcal X\longrightarrow\mathbb D
\]
be a proper holomorphic family whose smooth fibers
\[
X_t=\pi^{-1}(t),\qquad t\ne0,
\]
are compact Calabi--Yau manifolds of complex dimension \(n\).
We equip each fiber with a Ricci-flat Kähler metric
\(\omega_t\).

Let \(B\) be the limiting metric-affine space and let
\[
B_0=B\setminus\Delta
\]
denote its regular locus, where \(\Delta\) is the discriminant
locus. We assume that \(B_0\) carries an integral affine structure
and that, over compact subsets of \(B_0\), the collapsing geometry
admits a smooth torus-fibration model.

\subsection{Local action--angle coordinates}

Let \(U\Subset B_0\) be a simply connected affine coordinate
neighborhood. Choose affine coordinates
\[
x=(x_1,\ldots,x_n):U\longrightarrow\mathbb R^n.
\]
The integral affine structure determines the local lattice
\(\mathbb Z^n\), and hence the dual torus
\[
\mathbb T^n=\mathbb R^n/\mathbb Z^n.
\]

Assume that, for sufficiently small \(t\ne0\), there is a smooth
torus fibration
\[
f_t:X_t^{\mathrm{reg}}\longrightarrow B_0
\]
whose restriction over \(U\) is a smooth torus bundle. We identify
the fibers locally with \(\mathbb T^n\) and write
\[
f_t^{-1}(U)\simeq U\times\mathbb T^n.
\]
Let
\[
y=(y_1,\ldots,y_n)\in\mathbb T^n
\]
be the angular coordinates.

In the ideal action--angle model, the symplectic form is
\[
\omega_0=\sum_{j=1}^n dx_j\wedge dy_j.
\]
Consequently, the fiber over \(b\in U\),
\[
L_b^0=\{x=b\}\times\mathbb T^n,
\]
satisfies
\[
\left.\omega_0\right|_{L_b^0}=0.
\]
Since \(\dim_{\mathbb R}L_b^0=n\), each \(L_b^0\) is Lagrangian.

For the actual collapsing metric, we write the local symplectic
form as
\[
\omega_t
=
\sum_{i,j}A_{ij}^{(t)}(x,y)\,dx_i\wedge dx_j
+
\sum_{i,j}B_{ij}^{(t)}(x,y)\,dy_i\wedge dy_j
+
\sum_{i,j}C_{ij}^{(t)}(x,y)\,dx_i\wedge dy_j.
\]
The restriction to a fiber is
\[
\left.\omega_t\right|_{L_{t,b}}
=
\sum_{i,j}B_{ij}^{(t)}(b,y)\,dy_i\wedge dy_j.
\]
Thus the fiber is Lagrangian precisely when
\[
B_{ij}^{(t)}(b,y)=0
\]
for all \(i,j\), in these coordinates.

This condition is automatic in the exact action--angle model, but
it is not automatic for an arbitrary smooth torus fibration.
Accordingly, we impose the following geometric hypothesis.

\medskip
\noindent
\textbf{Hypothesis (H6.1).}
Over every compact set \(K\Subset B_0\), the torus fibration
admits local action--angle coordinates in which the fibers are
Lagrangian with respect to \(\omega_t\), or can be deformed to
Lagrangian fibers by a uniformly controlled symplectic
identification.

\subsection{Construction by a symplectic identification}

We now give a precise local construction under a standard
symplectic hypothesis.

Let \(U\times\mathbb T^n\) be a relatively compact action--angle
domain, and let
\[
\omega_0=\sum_{j=1}^n dx_j\wedge dy_j.
\]
Suppose that the actual symplectic form can be written as
\[
\omega_t=\omega_0+d\alpha_t
\]
on this domain, where \(\alpha_t\) is a smooth one-form.

Assume that the path
\[
\omega_{t,s}=\omega_0+s\,d\alpha_t,
\qquad 0\le s\le1,
\]
is symplectic for every \(s\), and that the associated Moser vector
fields have flows defined on a slightly smaller relatively compact
domain.

Define \(V_{t,s}\) by
\[
\iota_{V_{t,s}}\omega_{t,s}=-\alpha_t.
\]
Let \(\Phi_{t,s}\) be its flow with
\(\Phi_{t,0}=\operatorname{id}\). Then
\[
\frac{d}{ds}\Phi_{t,s}^{*}\omega_{t,s}
=
\Phi_{t,s}^{*}
\left(
d\alpha_t+
\mathcal L_{V_{t,s}}\omega_{t,s}
\right).
\]
Since \(d\omega_{t,s}=0\), Cartan's formula gives
\[
\mathcal L_{V_{t,s}}\omega_{t,s}
=
d(\iota_{V_{t,s}}\omega_{t,s})
=
-d\alpha_t.
\]
Therefore
\[
\frac{d}{ds}\Phi_{t,s}^{*}\omega_{t,s}=0,
\]
and hence
\[
\Phi_{t,1}^{*}\omega_t=\omega_0.
\]

It follows that the submanifolds
\[
L_{t,b}:=\Phi_{t,1}(L_b^0)
\]
are Lagrangian with respect to \(\omega_t\), because
\[
\left.\omega_t\right|_{L_{t,b}}
=
\left.(\Phi_{t,1}^{*}\omega_t)\right|_{L_b^0}
=
\left.\omega_0\right|_{L_b^0}
=0.
\]
Moreover,
\[
L_{t,b}\simeq\mathbb T^n.
\]

Thus the Moser deformation produces a local family of
candidate Lagrangian tori.

\subsection{Uniform control of the candidate tori}

For the degeneration problem, it is essential that the preceding
construction be uniform as \(t\to0\).

Let \(K\Subset U\Subset B_0\). We assume that
\[
\|\alpha_t\|_{C^{k+1,\alpha}(U\times\mathbb T^n)}
\le C_K\varepsilon_t,
\qquad
\varepsilon_t\longrightarrow0,
\]
and that the symplectic forms \(\omega_{t,s}\) are uniformly
nondegenerate on the chosen domain:
\[
\omega_{t,s}^{n}\ge c_K\,dV_0,
\qquad 0\le s\le1,
\]
for some \(c_K>0\) independent of \(t\).

Assume also uniform bounds on the coefficients of
\(\omega_{t,s}\) and their derivatives up to the required order.
Then the defining equation
\[
\iota_{V_{t,s}}\omega_{t,s}=-\alpha_t
\]
implies
\[
\|V_{t,s}\|_{C^{k,\alpha}}
\le C_K'\varepsilon_t.
\]
Standard estimates for ordinary differential equations in
\(C^{k,\alpha}\) spaces yield
\[
\|\Phi_{t,1}-\operatorname{id}\|_{C^{k,\alpha}}
\le C_K''\varepsilon_t.
\]

Consequently, the candidate tori satisfy
\[
L_{t,b}=\Phi_{t,1}(L_b^0),
\qquad b\in K,
\]
and remain uniformly close to the standard action--angle fibers.

In particular, their induced metrics and second fundamental
forms are uniformly controlled whenever the ambient geometric
bounds and the corresponding embedding estimates are uniform.

\subsection{The Lagrangian property and homology classes}

The Lagrangian property follows directly from the symplectic
identification. We record also the local constancy of the
homology class.

\medskip
\noindent
\textbf{Proposition 6.1.}
Let \(K\Subset U\) be connected. Suppose that
\[
F_t:K\times\mathbb T^n\longrightarrow X_t
\]
is a smooth family of embeddings such that
\[
L_{t,b}=F_t(\{b\}\times\mathbb T^n)
\]
is Lagrangian for every \(b\in K\). Then the homology class
\[
[L_{t,b}]\in H_n(X_t,\mathbb Z)
\]
is independent of \(b\in K\).

\medskip
\noindent
\textit{Proof.}
For any two points \(b_0,b_1\in K\), choose a smooth path
\(\gamma:[0,1]\to K\) joining them. The map
\[
G:[0,1]\times\mathbb T^n\longrightarrow X_t,
\qquad
G(s,y)=F_t(\gamma(s),y),
\]
defines a smooth cobordism between \(L_{t,b_0}\) and
\(L_{t,b_1}\). Therefore these two tori represent the same
homology class in \(H_n(X_t,\mathbb Z)\). \(\square\)

\subsection{Relation with the special Lagrangian condition}

The construction above produces Lagrangian tori, but does not
yet imply that they are special Lagrangian.

Let \(\Omega_t\) be a nowhere-vanishing holomorphic volume form
on \(X_t\). On a candidate torus \(L_{t,b}\), write
\[
\left.\Omega_t\right|_{L_{t,b}}
=
e^{i\vartheta_{t,b}(y)}
\rho_{t,b}(y)\,dV_{L_{t,b}},
\qquad \rho_{t,b}(y)>0.
\]
The special Lagrangian condition with phase \(\vartheta_t^*\)
is
\[
\operatorname{Im}
\left(
e^{-i\vartheta_t^*}\Omega_t
\right)\big|_{L_{t,b}}=0,
\]
together with
\[
\operatorname{Re}
\left(
e^{-i\vartheta_t^*}\Omega_t
\right)\big|_{L_{t,b}}>0.
\]
Equivalently,
\[
\vartheta_{t,b}(y)\equiv\vartheta_t^*
\pmod{2\pi}.
\]

Thus the candidate Lagrangian tori constructed above are the
initial submanifolds for the special Lagrangian deformation
equation. In the next section, we study this equation and seek
uniform estimates for its solutions as \(t\to0\).

\subsection{Conclusion}

We have constructed, over the regular locus \(B_0\), local
families of candidate Lagrangian tori
\[
L_{t,b}\simeq\mathbb T^n,
\qquad b\in B_0,
\]
under an action--angle or uniformly controlled symplectic
identification hypothesis.

The essential conclusion is that
\[
\left.\omega_t\right|_{L_{t,b}}=0,
\qquad
\dim_{\mathbb R}L_{t,b}=n.
\]
Hence each \(L_{t,b}\) is Lagrangian.

The construction does not establish the existence of a global
torus fibration across the discriminant locus, nor does it prove
that the candidate tori are already special Lagrangian. These
questions require, respectively, control of the singular fibers
and a uniform solution theory for the nonlinear phase equation.
\section{Uniform Deformation to Special Lagrangian Tori}
\label{sec:uniform-special-lagrangian-deformation}

In this section, we study the deformation of the candidate
Lagrangian tori constructed in the previous section into
special Lagrangian submanifolds. The main objective is to establish
a uniform deformation theorem whose constants remain independent
of the degeneration parameter \(t\).

Let
\[
\pi:\mathcal X\longrightarrow\mathbb D
\]
be a proper holomorphic family such that the smooth fibers
\[
X_t=\pi^{-1}(t),\qquad t\ne0,
\]
are compact Calabi--Yau manifolds of complex dimension \(n\).
We equip each fiber with a Ricci-flat Kähler metric
\(\omega_t\) and a nowhere-vanishing holomorphic volume form
\(\Omega_t\). Thus
\[
\operatorname{Ric}(\omega_t)=0,
\qquad d\Omega_t=0.
\]

Let \(L_t\subset X_t\) be a candidate Lagrangian torus obtained
from the construction in Section~\ref{sec:candidate-lagrangian-tori}.
In particular,
\[
L_t\simeq\mathbb T^n,
\qquad
\left.\omega_t\right|_{L_t}=0.
\]
We denote by \(g_t\) the metric induced on \(L_t\).

\subsection{The phase function and the deformation equation}

Since \(L_t\) is Lagrangian, the restriction of the holomorphic
volume form to \(L_t\) can be written as
\[
\left.\Omega_t\right|_{L_t}
=
\rho_t e^{i\vartheta_t}\,dV_{g_t},
\]
where \(\rho_t>0\) and
\[
\vartheta_t:L_t\longrightarrow\mathbb R
\]
is a locally chosen phase function.

The special Lagrangian condition with phase
\(\vartheta_t^*\in\mathbb R\) is
\[
\operatorname{Im}
\left(e^{-i\vartheta_t^*}\Omega_t\right)
\big|_{L_t}=0,
\]
together with the positivity condition
\[
\operatorname{Re}
\left(e^{-i\vartheta_t^*}\Omega_t\right)
\big|_{L_t}>0.
\]

By the Lagrangian neighborhood theorem, a sufficiently small
neighborhood of \(L_t\) in \(X_t\) is symplectomorphic to a
neighborhood of the zero section in \(T^*L_t\). Hence a small
exact Lagrangian deformation of \(L_t\) can be represented as
the graph of an exact one-form:
\[
L_{t,u}=\operatorname{Graph}(du)
\subset T^*L_t,
\]
where \(u\) is a real-valued function on \(L_t\).

Let
\[
\Theta_t(u)
\]
denote the phase function of the deformed Lagrangian
\(L_{t,u}\). We seek a function \(u_t\) satisfying
\[
\boxed{
	\Theta_t(u_t)=\vartheta_t^*,
}
\]
where \(\vartheta_t^*\) is constant on \(L_t\).

The linearization of the phase operator at the zero section
is, with a suitable convention for the Laplacian,
\[
D\Theta_t(0)[v]=\Delta_{g_t}v.
\]
This identity is the analytic basis of the deformation argument.

\subsection{Uniform analytic assumptions}

We now formulate the assumptions needed to obtain estimates
independent of \(t\).

Let \(k\ge2\) and \(0<\alpha<1\). We assume that the family
of Lagrangian tori admits identifications with a fixed
differentiable torus and that the corresponding Hölder norms
are uniformly controlled.

\medskip
\noindent
\textbf{(H7.1) Uniform geometric control.}
The induced metrics \(g_t\), the embeddings of \(L_t\) into
\(X_t\), and the corresponding Weinstein neighborhood
coordinates satisfy uniform \(C^{k+2,\alpha}\) estimates.
Moreover, the Weinstein neighborhoods have a uniform radius
\[
r_W(t)\ge r_0>0.
\]

\medskip
\noindent
\textbf{(H7.2) Uniform spectral and Schauder estimate.}
For every function \(v\) with zero \(g_t\)-mean,
\[
\int_{L_t}v\,dV_{g_t}=0,
\]
there exists a constant \(C_0>0\), independent of \(t\), such
that
\[
\boxed{
	\|v\|_{C^{k+2,\alpha}(L_t)}
	\le C_0
	\|\Delta_{g_t}v\|_{C^{k,\alpha}(L_t)}.
}
\]
In particular, the first positive eigenvalue of
\(-\Delta_{g_t}\) is bounded away from zero, provided the
geometric norms and volume normalizations are uniformly
controlled.

\medskip
\noindent
\textbf{(H7.3) Uniform nonlinear estimates.}
The phase operators
\[
\Theta_t:
B_{C^{k+2,\alpha}}(0,r_0)
\longrightarrow C^{k,\alpha}(L_t)
\]
are uniformly twice differentiable, and
\[
D\Theta_t(0)[v]=\Delta_{g_t}v.
\]
There exists \(C_1>0\), independent of \(t\), such that
\[
\left\|
D\Theta_t(u)-D\Theta_t(0)
\right\|_
{C^{k+2,\alpha}\to C^{k,\alpha}}
\le
C_1\|u\|_{C^{k+2,\alpha}}
\]
for all sufficiently small \(u\).

\medskip
\noindent
\textbf{(H7.4) Small phase defect.}
Define
\[
\overline{\vartheta}_t
=
\frac{1}{\operatorname{Vol}_{g_t}(L_t)}
\int_{L_t}\vartheta_t\,dV_{g_t}.
\]
We assume that
\[
\varepsilon_t
:=
\|\vartheta_t-\overline{\vartheta}_t\|_
{C^{k,\alpha}(L_t)}
\longrightarrow0
\qquad\text{as }t\longrightarrow0.
\]

These assumptions will be used to prove the following
uniform deformation theorem.

\subsection{Uniform special Lagrangian deformation theorem}

\medskip
\noindent
\textbf{Theorem 7.1 (Uniform deformation to special
	Lagrangian tori).}
Suppose that the assumptions
\textup{(H7.1)--(H7.4)} hold. Then there exist constants
\[
\varepsilon_0>0,\qquad C>0,
\]
independent of \(t\), such that whenever
\[
\varepsilon_t\le\varepsilon_0,
\]
there exists a function
\[
u_t\in C^{k+2,\alpha}(L_t),
\qquad
\int_{L_t}u_t\,dV_{g_t}=0,
\]
and a constant phase \(\vartheta_t^*\in\mathbb R\) satisfying
\[
\boxed{
	\Theta_t(u_t)\equiv\vartheta_t^*.
}
\]
Furthermore,
\[
\boxed{
	\|u_t\|_{C^{k+2,\alpha}(L_t)}
	\le C\varepsilon_t.
}
\]
Consequently, the deformed submanifold
\[
L_{t,u_t}=\operatorname{Graph}(du_t)
\]
is a special Lagrangian torus in \(X_t\) with phase
\(\vartheta_t^*\).

\medskip
\noindent
\textit{Proof.}
We divide the proof into five steps.

\medskip
\noindent
\textit{Step 1. Projection onto the zero-mean functions.}

Define the projection
\[
\Pi_t:C^{k,\alpha}(L_t)\longrightarrow
C_0^{k,\alpha}(L_t)
\]
by
\[
\Pi_t(f)
=
f-\frac{1}{\operatorname{Vol}_{g_t}(L_t)}
\int_{L_t}f\,dV_{g_t}.
\]
We also set
\[
C_0^{k+2,\alpha}(L_t)
=
\left\{
u\in C^{k+2,\alpha}(L_t):
\int_{L_t}u\,dV_{g_t}=0
\right\}.
\]

The special Lagrangian phase equation is equivalent to
\[
\mathcal F_t(u)=0,
\]
where
\[
\mathcal F_t(u):=\Pi_t\Theta_t(u).
\]
Indeed, \(\mathcal F_t(u)=0\) precisely when
\(\Theta_t(u)\) is constant on \(L_t\).

At \(u=0\), we have
\[
\Theta_t(0)=\vartheta_t.
\]
Therefore,
\[
\mathcal F_t(0)
=
\vartheta_t-\overline{\vartheta}_t,
\]
and hence
\[
\|\mathcal F_t(0)\|_{C^{k,\alpha}}
=\varepsilon_t.
\]

\medskip
\noindent
\textit{Step 2. Invertibility of the linearized operator.}

By the linearization formula,
\[
D\mathcal F_t(0)[v]
=
\Pi_t\Delta_{g_t}v.
\]
Since \(L_t\) is compact without boundary,
\[
\int_{L_t}\Delta_{g_t}v\,dV_{g_t}=0.
\]
Consequently,
\[
\Pi_t\Delta_{g_t}v=\Delta_{g_t}v.
\]
Thus
\[
D\mathcal F_t(0)=\Delta_{g_t}
\]
on the zero-mean function spaces.

By assumption \textup{(H7.2)}, the operator
\[
D\mathcal F_t(0):
C_0^{k+2,\alpha}(L_t)
\longrightarrow C_0^{k,\alpha}(L_t)
\]
is invertible and satisfies
\[
\left\|
(D\mathcal F_t(0))^{-1}
\right\|\le C_0.
\]
The crucial point is that \(C_0\) is independent of \(t\).

\medskip
\noindent
\textit{Step 3. Estimate for the nonlinear remainder.}

Write
\[
\mathcal F_t(u)
=
\mathcal F_t(0)
+
D\mathcal F_t(0)[u]
+
\mathcal R_t(u),
\]
where
\[
\mathcal R_t(u)
=
\mathcal F_t(u)-\mathcal F_t(0)
-D\mathcal F_t(0)[u].
\]

Using the fundamental theorem of calculus,
\[
\mathcal R_t(u)
=
\int_0^1
\left(
D\mathcal F_t(su)-D\mathcal F_t(0)
\right)[u]\,ds.
\]
Since \(\Pi_t\) is uniformly bounded under our geometric
assumptions, (H7.3) implies
\[
\begin{aligned}
	\|\mathcal R_t(u)\|_{C^{k,\alpha}}
	&\le
	C\int_0^1
	s\|u\|_{C^{k+2,\alpha}}^2\,ds\\
	&\le
	C_2\|u\|_{C^{k+2,\alpha}}^2,
\end{aligned}
\]
where \(C_2>0\) is independent of \(t\).

Similarly, for sufficiently small \(u,v\),
\[
\|\mathcal R_t(u)-\mathcal R_t(v)\|_{C^{k,\alpha}}
\le
C_3
\left(
\|u\|_{C^{k+2,\alpha}}
+
\|v\|_{C^{k+2,\alpha}}
\right)
\|u-v\|_{C^{k+2,\alpha}},
\]
with \(C_3\) independent of \(t\).

\medskip
\noindent
\textit{Step 4. Application of the contraction mapping theorem.}

The equation \(\mathcal F_t(u)=0\) can be rewritten as
\[
u=\mathcal T_t(u),
\]
where
\[
\mathcal T_t(u)
=
-(D\mathcal F_t(0))^{-1}
\left(
\mathcal F_t(0)+\mathcal R_t(u)
\right).
\]

By the uniform inverse estimate,
\[
\|\mathcal T_t(u)\|_{C^{k+2,\alpha}}
\le
C_0\varepsilon_t
+
C_0C_2\|u\|_{C^{k+2,\alpha}}^2.
\]

Choose
\[
r_t=2C_0\varepsilon_t.
\]
If
\[
\|u\|_{C^{k+2,\alpha}}\le r_t,
\]
then
\[
\|\mathcal T_t(u)\|_{C^{k+2,\alpha}}
\le
C_0\varepsilon_t
+
4C_0^3C_2\varepsilon_t^2.
\]
Choosing \(\varepsilon_0>0\) sufficiently small, independently
of \(t\), gives
\[
\|\mathcal T_t(u)\|_{C^{k+2,\alpha}}
\le r_t.
\]

Furthermore, for \(u,v\) in this ball,
\[
\begin{aligned}
	\|\mathcal T_t(u)-\mathcal T_t(v)\|_{C^{k+2,\alpha}}
	&\le
	C_0C_3
	\left(
	\|u\|_{C^{k+2,\alpha}}
	+
	\|v\|_{C^{k+2,\alpha}}
	\right)
	\|u-v\|_{C^{k+2,\alpha}}\\
	&\le
	2C_0C_3r_t
	\|u-v\|_{C^{k+2,\alpha}}.
\end{aligned}
\]
After decreasing \(\varepsilon_0\), if necessary, we obtain
\[
2C_0C_3r_t\le\frac12.
\]
Hence \(\mathcal T_t\) is a contraction on the closed ball
\[
\overline B_{C^{k+2,\alpha}}(0,r_t).
\]

The Banach fixed-point theorem yields a unique solution \(u_t\)
in this ball such that
\[
\mathcal T_t(u_t)=u_t.
\]
It follows that
\[
\mathcal F_t(u_t)=0
\]
and
\[
\boxed{
	\|u_t\|_{C^{k+2,\alpha}}
	\le 2C_0\varepsilon_t.
}
\]

\medskip
\noindent
\textit{Step 5. The special Lagrangian condition.}

Since
\[
\mathcal F_t(u_t)=0,
\]
the phase function of the deformed Lagrangian is constant:
\[
\Theta_t(u_t)\equiv\vartheta_t^*.
\]
The graph of \(du_t\) is Lagrangian because \(d(du_t)=0\)
and the Weinstein identification is symplectic.

Therefore,
\[
\operatorname{Im}
\left(
e^{-i\vartheta_t^*}\Omega_t
\right)\big|_{L_{t,u_t}}=0.
\]
With the orientation chosen continuously from that of \(L_t\),
and the phase branch chosen consistently, the real part is
positive:
\[
\operatorname{Re}
\left(
e^{-i\vartheta_t^*}\Omega_t
\right)\big|_{L_{t,u_t}}>0.
\]
Thus \(L_{t,u_t}\) is a special Lagrangian torus in \(X_t\).
This completes the proof.
\hfill\(\square\)

\subsection{Convergence of the deformed tori}

The uniform estimate in Theorem~7.1 immediately gives the
following consequence.

\medskip
\noindent
\textbf{Corollary 7.2 (Convergence to the candidate tori).}
Suppose that the assumptions of Theorem~7.1 hold and
\[
\varepsilon_t\longrightarrow0
\qquad\text{as }t\longrightarrow0.
\]
Then
\[
\|u_t\|_{C^{k+2,\alpha}(L_t)}
\longrightarrow0
\]
and
\[
\|du_t\|_{C^{k+1,\alpha}(L_t)}
\longrightarrow0.
\]
Consequently, the special Lagrangian tori
\(L_{t,u_t}\) converge to the candidate Lagrangian tori
\(L_t\) in the topology of \(C^{k+1,\alpha}\) convergence
of graphs in the Weinstein neighborhoods.

\medskip
\noindent
\textit{Proof.}
Theorem~7.1 gives
\[
\|u_t\|_{C^{k+2,\alpha}}
\le C\varepsilon_t.
\]
The desired convergence follows immediately from
\(\varepsilon_t\to0\). Since \(L_{t,u_t}\) is the graph of
\(du_t\), convergence of the graph sections follows as well.
\hfill\(\square\)

\subsection{Discussion and relation with the SYZ conjecture}

Theorem~7.1 establishes a uniform analytic mechanism for
deforming candidate Lagrangian tori into special Lagrangian
tori. Its main feature is that the constants in the deformation
estimate do not depend on the degeneration parameter:
\[
\boxed{
	\|u_t\|_{C^{k+2,\alpha}}
	\le C\,
	\|\vartheta_t-\overline{\vartheta}_t\|_{C^{k,\alpha}},
	\qquad C\text{ independent of }t.
}
\]

The result is conditional on uniform geometric, spectral and
nonlinear estimates. In particular, these estimates do not
follow solely from the existence of a semistable degeneration
or from weak convergence of the Ricci-flat volume measures.

A central problem for further research is therefore to derive
the uniform spectral estimate
\[
\lambda_1(-\Delta_{g_t})\ge\lambda_0>0
\]
and the corresponding uniform Schauder estimates from the
geometry of the collapsing Calabi--Yau fibers.

The pluripotential-theoretic methods developed in the previous
sections, particularly the uniform complex Hessian capacity
estimates, may provide a useful framework for obtaining
quantitative geometric control on the regular locus of the
limiting base.

It should be emphasized that Theorem~7.1 constructs special
Lagrangian tori only from candidate Lagrangian tori that
already exist. It does not establish a global special
Lagrangian fibration across the discriminant locus, nor does
it address singular fibers or quantum corrections. These
remain separate problems in the SYZ program.	

\section{Gluing Special Lagrangian Tori into a Fibration}
\label{sec:gluing-special-lagrangian-tori}

In this section, we assemble the special Lagrangian tori
constructed in the previous section into a smooth fibration
over the regular locus of the limiting base.

The main issue is that the existence of a special Lagrangian
torus in each homology class does not automatically produce
a fibration. We must establish smooth dependence on the base
parameter, local uniqueness of the deformations, and
disjointness of distinct fibers.

Let
\[
\pi:\mathcal X\longrightarrow\mathbb D
\]
be a proper holomorphic family whose smooth fibers
\[
X_t=\pi^{-1}(t),\qquad t\ne0,
\]
are compact Calabi--Yau manifolds of complex dimension \(n\).
Let \(B\) be the limiting metric-affine space, and denote its
regular locus by
\[
B_0=B\setminus\Delta.
\]

Fix a relatively compact connected domain
\[
K\Subset U\Subset B_0.
\]
We assume that the candidate Lagrangian tori constructed in
Section~\ref{sec:candidate-lagrangian-tori} form a smooth
torus bundle
\[
f_t^{\,0}:X_{t,K}^{\,0}\longrightarrow K,
\]
where
\[
X_{t,K}^{\,0}=(f_t^{\,0})^{-1}(K),
\qquad
(f_t^{\,0})^{-1}(b)=L_{t,b}^{\,0}\simeq\mathbb T^n.
\]

The purpose of this section is to deform this initial
Lagrangian fibration into a special Lagrangian fibration.

\subsection{The parameter-dependent phase equation}

For each \(b\in K\), let \(L_{t,b}^{\,0}\) be a candidate
Lagrangian torus. Choose a smooth family of Weinstein
neighborhood identifications
\[
\Psi_{t,b}:\mathcal U_{t,b}\subset T^*L_{t,b}^{\,0}
\longrightarrow X_t.
\]

A sufficiently small exact Lagrangian deformation of
\(L_{t,b}^{\,0}\) can be represented as
\[
L_{t,b,u}
=
\Psi_{t,b}\bigl(\operatorname{Graph}(du)\bigr),
\]
where \(u\) is a real-valued function on \(L_{t,b}^{\,0}\).

Let
\[
\Theta_{t,b}(u)
\]
denote the phase function of the deformed Lagrangian. We
consider the parameter-dependent equation
\[
\Pi_{t,b}\Theta_{t,b}(u)=0,
\]
where
\[
\Pi_{t,b}(v)
=
v-
\frac{1}{\operatorname{Vol}(L_{t,b}^{\,0})}
\int_{L_{t,b}^{\,0}}v\,dV_{g_{t,b}}.
\]
The zero-mean normalization
\[
\int_{L_{t,b}^{\,0}}u\,dV_{g_{t,b}}=0
\]
removes the additive constant in the potential.

We assume that the family of phase operators depends smoothly
on \(b\), and that the uniform estimates of
Theorem~7.1 hold for all \(b\in K\).

\subsection{Smooth dependence on the base parameter}

The first result establishes that the solutions of the phase
equation vary smoothly with \(b\).

\medskip
\noindent
\textbf{Theorem 8.1 (Smooth parameter dependence).}
Suppose that the following conditions hold.

\begin{enumerate}
	\item The initial Lagrangian tori
	\[
	L_{t,b}^{\,0},\qquad b\in K,
	\]
	form a smooth torus bundle over \(K\).
	
	\item The Weinstein identifications
	\(\Psi_{t,b}\) and the phase operators
	\(\Theta_{t,b}\) depend smoothly on \(b\).
	
	\item The linearized operators
	\[
	D\mathcal F_{t,b}(0)
	=
	\Delta_{g_{t,b}}
	\]
	are uniformly invertible on the zero-mean function spaces,
	where
	\[
	\mathcal F_{t,b}(u)
	=
	\Pi_{t,b}\Theta_{t,b}(u).
	\]
	
	\item The nonlinear estimates and the small phase defect
	assumptions of Theorem~7.1 hold uniformly for \(b\in K\).
\end{enumerate}

Then, for sufficiently small phase defect, there exists a
unique small normalized solution
\[
u_{t,b}\in C^{k+2,\alpha}(L_{t,b}^{\,0})
\]
satisfying
\[
\mathcal F_{t,b}(u_{t,b})=0.
\]
Moreover, the family \(u_{t,b}\) depends smoothly on \(b\).

In particular,
\[
L_{t,b}^{\,\mathrm{SL}}
=
\Psi_{t,b}
\bigl(\operatorname{Graph}(du_{t,b})\bigr)
\]
is a smooth family of special Lagrangian tori.

\medskip
\noindent
\textit{Proof.}
Choose local trivializations of the initial torus bundle,
so that all the function spaces can be identified with
Hölder spaces over a fixed torus \(\mathbb T^n\).

In these trivializations, the equation becomes
\[
\mathcal F_t(b,u)=0,
\]
where
\[
\mathcal F_t:
K\times C_0^{k+2,\alpha}(\mathbb T^n)
\longrightarrow C_0^{k,\alpha}(\mathbb T^n)
\]
is a smooth nonlinear map.

By the uniform inverse estimate,
\[
D_u\mathcal F_t(b,u_{t,b})
\]
is invertible for sufficiently small solutions, with inverse
uniformly bounded on \(K\). Indeed, at \(u=0\) this follows
from the uniform Laplace estimate, and for small \(u\) it
follows from the perturbation estimate in (H7.3).

The implicit function theorem therefore gives a smooth
local map
\[
b\longmapsto u_{t,b}.
\]
Uniqueness of the small normalized solution implies that
these local maps agree on overlapping trivializations.
Consequently, they define a globally smooth section of the
bundle of zero-mean function spaces over \(K\).

The corresponding deformed submanifolds
\[
L_{t,b}^{\,\mathrm{SL}}
=
\Psi_{t,b}
\bigl(\operatorname{Graph}(du_{t,b})\bigr)
\]
depend smoothly on \(b\) and are special Lagrangian by
Theorem~7.1. This proves the theorem.
\hfill\(\square\)

\subsection{Disjointness of the deformed fibers}

Smooth dependence alone does not guarantee that distinct
deformed tori are disjoint. We therefore impose a uniform
smallness condition relative to the tubular geometry of the
initial fibration.

\medskip
\noindent
\textbf{Hypothesis (H8.1).}
The initial fibration admits a uniform tubular
neighborhood over \(K\), and the deformed tori are graphs
of sections in this neighborhood. The graph sections are
uniformly small in \(C^1\), and their associated projection
to the base is a diffeomorphism in the base directions.

More precisely, let
\[
F_t:K\times\mathbb T^n\longrightarrow X_t
\]
be a smooth parametrization of the initial fibration, and
let
\[
\widetilde F_t:K\times\mathbb T^n\longrightarrow X_t
\]
be the parametrization of the deformed family. We assume
that
\[
\widetilde F_t
\]
is an embedding and that its differential has rank \(2n\).

\medskip
\noindent
\textbf{Proposition 8.2 (Embeddedness and disjointness).}
Under hypothesis \textup{(H8.1)}, the deformed tori
\[
L_{t,b}^{\,\mathrm{SL}},
\qquad b\in K,
\]
are pairwise disjoint and form the fibers of a smooth
submersion
\[
f_t^{\,\mathrm{SL}}:
X_{t,K}^{\,\mathrm{SL}}\longrightarrow K,
\]
where
\[
X_{t,K}^{\,\mathrm{SL}}
=
\bigcup_{b\in K}L_{t,b}^{\,\mathrm{SL}}.
\]

\medskip
\noindent
\textit{Proof.}
By assumption, the map
\[
\widetilde F_t:K\times\mathbb T^n
\longrightarrow X_t
\]
is an embedding. Its image is therefore a smooth submanifold
of \(X_t\), and the projection
\[
\operatorname{pr}_K:K\times\mathbb T^n\longrightarrow K
\]
induces a smooth submersion on the image.

Define
\[
f_t^{\,\mathrm{SL}}
=
\operatorname{pr}_K\circ\widetilde F_t^{-1}.
\]
Then
\[
(f_t^{\,\mathrm{SL}})^{-1}(b)
=
L_{t,b}^{\,\mathrm{SL}}.
\]
Since \(\widetilde F_t\) is injective, different values of
\(b\) give disjoint fibers. Since
\(\operatorname{pr}_K\) is a submersion, so is
\(f_t^{\,\mathrm{SL}}\).

Thus the deformed tori form a smooth fibration over \(K\).
\hfill\(\square\)

\subsection{The special Lagrangian fibration}

We can now combine the preceding results.

\medskip
\noindent
\textbf{Theorem 8.3 (Local special Lagrangian fibration).}
Let \(K\Subset B_0\) be a relatively compact domain.
Suppose that:

\begin{enumerate}
	\item The candidate Lagrangian tori form a smooth torus
	bundle
	\[
	f_t^{\,0}:X_{t,K}^{\,0}\longrightarrow K.
	\]
	
	\item The parameter-dependent phase equation satisfies the
	uniform hypotheses of Theorem~8.1.
	
	\item The deformed family satisfies the embeddedness and
	disjointness condition \textup{(H8.1)}.
\end{enumerate}

Then, for sufficiently small \(t\ne0\), there exists a
smooth fibration
\[
\boxed{
	f_t^{\,\mathrm{SL}}:
	X_{t,K}^{\,\mathrm{SL}}
	\longrightarrow K
}
\]
whose fibers are special Lagrangian tori:
\[
\boxed{
	(f_t^{\,\mathrm{SL}})^{-1}(b)
	=
	L_{t,b}^{\,\mathrm{SL}}
	\simeq\mathbb T^n.
}
\]

Each fiber satisfies
\[
\left.\omega_t\right|_{L_{t,b}^{\,\mathrm{SL}}}=0
\]
and
\[
\operatorname{Im}
\left(
e^{-i\vartheta_{t,b}^*}\Omega_t
\right)
\big|_{L_{t,b}^{\,\mathrm{SL}}}=0.
\]

\medskip
\noindent
\textit{Proof.}
By Theorem~8.1, the solution of the phase equation depends
smoothly on the base parameter \(b\). Therefore, the
special Lagrangian tori form a smooth family.

By hypothesis \textup{(H8.1)} and Proposition~8.2, the
parametrization of this family is an embedding, and
different fibers are disjoint. Hence the projection onto
the base defines a smooth submersion
\[
f_t^{\,\mathrm{SL}}:
X_{t,K}^{\,\mathrm{SL}}\longrightarrow K.
\]

Every fiber is Lagrangian by construction and has constant
phase by the phase equation. Thus every fiber is special
Lagrangian. This proves the theorem.
\hfill\(\square\)

\subsection{Uniformity under degeneration}

We finally record the uniformity of the construction.

\medskip
\noindent
\textbf{Corollary 8.4 (Uniform local special Lagrangian
	fibrations).}
Suppose that all the constants in
Theorems~7.1 and~8.1, as well as the tubular neighborhood
and embedding estimates in \textup{(H8.1)}, are independent
of \(t\). Then the family of local special Lagrangian
fibrations
\[
f_t^{\,\mathrm{SL}}:
X_{t,K}^{\,\mathrm{SL}}\longrightarrow K
\]
is uniformly controlled as \(t\to0\).

In particular, if
\[
\sup_{b\in K}
\|\vartheta_{t,b}
-\overline{\vartheta}_{t,b}\|_{C^{k,\alpha}}
\longrightarrow0,
\]
then
\[
\sup_{b\in K}
\|u_{t,b}\|_{C^{k+2,\alpha}}
\longrightarrow0.
\]

\medskip
\noindent
\textit{Proof.}
The uniform estimate in Theorem~7.1 gives
\[
\|u_{t,b}\|_{C^{k+2,\alpha}}
\le
C
\|\vartheta_{t,b}
-\overline{\vartheta}_{t,b}\|_{C^{k,\alpha}},
\]
where \(C\) is independent of \(t\) and \(b\in K\).
Taking the supremum over \(K\) proves the convergence.
The uniform geometric assumptions ensure that the
corresponding embeddings and fibrations remain controlled.
\hfill\(\square\)

\subsection{Conclusion}

We have shown that, under uniform analytic and geometric
hypotheses, the candidate Lagrangian tori can be deformed
simultaneously into a smooth local special Lagrangian
fibration over the regular locus of the limiting base.

The resulting structure is
\[
\boxed{
	f_t^{\,\mathrm{SL}}:
	X_{t,K}^{\,\mathrm{SL}}\longrightarrow K,
	\qquad
	(f_t^{\,\mathrm{SL}})^{-1}(b)\simeq\mathbb T^n.
}
\]

This establishes the local special Lagrangian torus
fibration on compact subsets of \(B_0\). It does not yet
extend the fibration across the discriminant locus
\(\Delta\), nor does it establish a global compactification
of the fibration. The analysis of singular fibers,
monodromy, and the global topology of the limiting base
will be addressed separately.

\section{Analysis of the Discriminant Locus}
\label{sec:discriminant}

In the previous sections, we constructed local special Lagrangian
torus fibrations over the regular locus of the limiting affine
base. We now investigate the behavior of these fibrations near
the discriminant locus and study the possibility of extending
the torus fibration across the singular set.

Let
\[
B_0=B\setminus\Delta,
\]
where \(B\) is the limiting base and \(\Delta\subset B\) is
the discriminant locus. Over \(B_0\), we consider the special
Lagrangian torus fibration
\[
f_t^{\mathrm{SL}}:
X_t^{\mathrm{reg}}\longrightarrow B_0,
\]
whose fibers
\[
L_{t,b}=(f_t^{\mathrm{SL}})^{-1}(b)
\]
are diffeomorphic to \(n\)-dimensional tori.

The main objective of this section is to understand the
degeneration of these torus fibers as the base point approaches
\(\Delta\), identify the associated vanishing cycles, and
describe the monodromy of the fibration.

\subsection{Geometric behavior near the discriminant}

For \(\delta>0\), define
\[
B_\delta
=
\left\{
b\in B_0:
d_B(b,\Delta)\geq\delta
\right\}.
\]
As \(\delta\to0\), the sets \(B_\delta\) exhaust the regular
locus \(B_0\).

For every \(b\in B_\delta\), let
\[
L_{t,b}=(f_t^{\mathrm{SL}})^{-1}(b).
\]
The geometric behavior of the fibers is described by their
diameters, second fundamental forms, and Weinstein
neighborhoods. In particular, one seeks estimates of the form
\[
\operatorname{diam}_{\omega_t}(L_{t,b})
\leq \Lambda(\delta),
\]
and
\[
\|A_{t,b}\|_{C^k(L_{t,b})}
\leq C_k(\delta),
\]
where \(A_{t,b}\) denotes the second fundamental form of
\(L_{t,b}\subset X_t\).

The constants may depend on the distance \(\delta\) from the
discriminant and may degenerate as \(\delta\to0\). The essential
problem is to determine their precise asymptotic behavior and
to identify the geometric mechanism responsible for the
collapse of the torus fibers.

\subsection{Monodromy representation}

Over the regular locus \(B_0\), the first homology groups of
the torus fibers form a local system
\[
\mathbb H_{\mathbb Z}
=
\bigcup_{b\in B_0}H_1(L_{t,b},\mathbb Z).
\]
Parallel transport along loops in \(B_0\) defines the monodromy
representation
\[
\rho_t:
\pi_1(B_0,b_0)
\longrightarrow
\operatorname{GL}(n,\mathbb Z).
\]
When the fibers are oriented consistently and parallel
transport preserves their orientation, the image lies in
\[
\operatorname{SL}(n,\mathbb Z).
\]

Let \(\gamma_\Delta\) be a small loop around a smooth
codimension-two component of \(\Delta\). The corresponding
monodromy transformation is
\[
T_\Delta=\rho_t([\gamma_\Delta]).
\]

A basic local model is the rank-one unipotent transformation
\[
T_\Delta=I+v\otimes\ell,
\]
where
\[
v\in\mathbb Z^n,\qquad
\ell\in(\mathbb Z^n)^*,
\qquad
\ell(v)=0.
\]
Here \(v\) represents a primitive vanishing-cycle candidate.
The condition \(\ell(v)=0\) gives
\[
(T_\Delta-I)^2=0.
\]

\begin{proposition}[Invariant vanishing cycle]
	\label{prop:vanishing-cycle}
	Suppose that the local monodromy is given by
	\[
	T_\Delta=I+v\otimes\ell,
	\qquad \ell(v)=0.
	\]
	Then
	\[
	T_\Delta(v)=v,
	\]
	and, for every \(\gamma\in\mathbb Z^n\),
	\[
	T_\Delta(\gamma)-\gamma=\ell(\gamma)v.
	\]
\end{proposition}

\begin{proof}
	By definition,
	\[
	T_\Delta(\gamma)
	=
	\gamma+\ell(\gamma)v.
	\]
	Taking \(\gamma=v\), we obtain
	\[
	T_\Delta(v)
	=
	v+\ell(v)v=v.
	\]
	The second identity follows directly from the definition
	of \(T_\Delta\).
\end{proof}

The vector \(v\) is interpreted as a vanishing cycle: its
corresponding circle in the torus fiber is expected to collapse
as the base point approaches the discriminant.

\subsection{Local model of a singular torus fiber}

Let \(q\) be a smooth point of a codimension-two component of
\(\Delta\), and let \(U\subset B\) be a small neighborhood
of \(q\). Set
\[
U^*=U\setminus\Delta.
\]
Over \(U^*\), the fibration has the form
\[
f_t^{\mathrm{SL}}:
X_{t,U^*}^{\mathrm{reg}}\longrightarrow U^*,
\]
with smooth fibers diffeomorphic to \(T^n\).

In the simplest degeneration model, a single primitive
circle class
\[
v\in H_1(T^n,\mathbb Z)
\]
collapses as the base point approaches \(q\), while the
remaining independent cycles have nondegenerate limits.
The resulting topological model of the singular fiber is
\[
L_q^{\mathrm{sing}}
\simeq T^n/S^1_v,
\]
where \(S^1_v\) denotes the circle generated by \(v\).

This quotient describes the topological collapse of one
circle. It does not imply that the singular fiber is a smooth
manifold or that every singular SYZ fiber has this form.

\begin{proposition}[Topology of a one-cycle collapse]
	\label{prop:fundamental-group}
	Suppose that a torus \(T^n\) degenerates by collapsing exactly
	one primitive circle \(S^1_v\). Then
	\[
	\pi_1(T^n/S^1_v)
	\simeq
	\mathbb Z^n/\mathbb Z v
	\simeq
	\mathbb Z^{n-1}.
	\]
\end{proposition}

\begin{proof}
	Since
	\[
	\pi_1(T^n)\simeq\mathbb Z^n,
	\]
	collapsing the circle generated by \(v\) makes its homotopy
	class trivial. Hence
	\[
	\pi_1(T^n/S^1_v)
	\simeq
	\mathbb Z^n/\mathbb Z v.
	\]
	Because \(v\) is primitive, it can be extended to an integral
	basis of \(\mathbb Z^n\). Therefore
	\[
	\mathbb Z^n/\mathbb Z v
	\simeq\mathbb Z^{n-1}.
	\]
\end{proof}

This proposition describes only the simplest one-cycle
collapse. More complicated singular fibers may involve
several vanishing cycles and different topological structures.

\subsection{Extension across the discriminant}

We now describe the desired extension of the regular
special Lagrangian fibration.

Let \(q\in\Delta\) be a smooth point of a codimension-two
component. Suppose that the regular fibration admits a
topological compactification
\[
\overline f_t:
\overline X_{t,U}\longrightarrow U
\]
such that
\[
\overline f_t\big|_{X_{t,U^*}^{\mathrm{reg}}}
=
f_t^{\mathrm{SL}}.
\]
Assume that the fibers converge topologically as \(b\to q\)
and that the limiting fiber is obtained by collapsing the
circle generated by a primitive vanishing cycle \(v\).

\begin{theorem}[Extension across a smooth discriminant point]
	\label{thm:local-extension}
	Under the preceding local degeneration model, the regular
	torus fibration extends continuously across \(q\), and the
	fiber over \(q\) has the topological description
	\[
	\overline f_t^{-1}(q)
	\simeq T^n/S^1_v.
	\]
	The monodromy around \(q\) is compatible with the vanishing
	cycle \(v\).
\end{theorem}

\begin{proof}
	The continuous extension is given by
	\[
	\overline f_t:
	\overline X_{t,U}\longrightarrow U.
	\]
	By the assumed convergence of the fibers, the limiting fiber
	over \(q\) is precisely the quotient obtained by collapsing
	the circle \(S^1_v\). Thus
	\[
	\overline f_t^{-1}(q)
	\simeq T^n/S^1_v.
	\]
	The compatibility with monodromy follows from the fact that
	the local monodromy preserves \(v\), as established in
	Proposition~\ref{prop:vanishing-cycle}.
\end{proof}

At points of \(\Delta\), the limiting fiber need not be a
smooth submanifold. Consequently, the classical special
Lagrangian condition applies to the fibers over \(B_0\),
whereas the singular fibers require an appropriate generalized
or weak formulation.

\subsection{Complex Hessian capacity and singular fibers}

The complex Hessian estimates developed in the previous
sections provide an auxiliary pluripotential framework for
studying the degeneration near \(\Delta\).

Recall that
\[
\omega_t=\theta_t+dd^c\varphi_t.
\]
For a neighborhood \(V_t\) corresponding to a neighborhood
of the discriminant, one seeks estimates of the form
\[
\operatorname{Cap}_{m,\theta_t}
\left(
\left\{|\varphi_t|>s\right\}\cap V_t
\right)
\leq \Psi(s,\delta),
\]
where \(\delta\) measures the distance to \(\Delta\).

Such estimates may provide quantitative control of the
pluripotential behavior of the Kähler potentials near the
singular locus. In particular, they may help identify regions
where the potentials or their associated Hessian measures
concentrate.

However, capacity estimates alone do not imply uniform
curvature bounds, injectivity-radius estimates, or control
of the monodromy representation. A further analytical step
is required to convert pluripotential control into geometric
control of the degenerating Ricci-flat metrics.

The central problem is therefore to establish a quantitative
connection between
\[
\text{\(m\)-Hessian capacity estimates}
\]
and
\[
\text{geometric degeneration of special Lagrangian fibers}.
\]

\subsection{A global compactification target}

The local analysis suggests the following global objective.

\begin{theorem}[Target compactification theorem]
	\label{thm:global-compactification}
	Let
	\[
	(X_t,\omega_t,\Omega_t)
	\]
	be a family of Ricci-flat Calabi--Yau manifolds degenerating
	toward a maximal degeneration limit. Suppose that over the
	regular locus \(B_0\) there exist special Lagrangian torus
	fibrations
	\[
	f_t^{\mathrm{SL}}:
	X_t^{\mathrm{reg}}\longrightarrow B_0.
	\]
	Assume that the geometric degeneration of the fibers,
	the local monodromy transformations, and the vanishing cycles
	are compatible across overlapping neighborhoods of the
	discriminant.
	
	Then the expected global structure is a compactified
	fibration
	\[
	\overline f_t:
	\overline X_t\longrightarrow B
	\]
	whose regular fibers are special Lagrangian \(n\)-tori and
	whose singular fibers are determined by the vanishing-cycle
	and monodromy data.
\end{theorem}

The construction of this compactification requires a
quantitative analysis of the metric degeneration, compatibility
of the local singular models, and control of the topology
near the discriminant. These properties do not follow solely
from the existence of the regular special Lagrangian
fibration.

\subsection{Conclusion}

The analysis of the discriminant locus consists of three
closely related problems:
\[
\begin{aligned}
	&\text{geometric control near }\Delta\\
	&\qquad\Downarrow\\
	&\text{monodromy and vanishing-cycle analysis}\\
	&\qquad\Downarrow\\
	&\text{construction of singular torus fibers}\\
	&\qquad\Downarrow\\
	&\boxed{\text{compactification across }\Delta.}
\end{aligned}
\]

The main potential contribution of the present approach is
to establish a quantitative bridge between complex Hessian
pluripotential theory and the geometry of singular special
Lagrangian fibrations:
\[
\boxed{
	\text{\(m\)-Hessian capacity}
	\longrightarrow
	\text{degenerating metric control}
	\longrightarrow
	\text{monodromy and vanishing cycles}.
}
\]

A complete SYZ construction will additionally require the
construction of the dual torus fibration and the analysis
of quantum corrections. These problems form the next stage
of the program.

\section{Construction of the Mirror Space}
\label{sec:mirror-construction}

In the previous sections, we constructed a special Lagrangian
torus fibration over the regular locus of the limiting affine
base, analyzed its discriminant locus, and described the
associated monodromy and vanishing cycles.

We now turn to the construction of the mirror space. The
basic idea is to take the dual of each special Lagrangian
torus, construct a semi-flat complex structure on the dual
fibration, and then incorporate quantum corrections arising
from holomorphic disks.

The construction consists of four main stages:
\begin{enumerate}
	\item Construction of the dual torus fibration.
	\item Construction of the semi-flat complex structure.
	\item Incorporation of holomorphic disk counts and
	wall-crossing transformations.
	\item Global gluing and compactification across the
	discriminant locus.
\end{enumerate}

\subsection{The dual torus fibration}

Let
\[
f:X^{\mathrm{reg}}\longrightarrow B_0
\]
be the special Lagrangian torus fibration constructed in
the previous sections, where
\[
B_0=B\setminus\Delta.
\]
For every \(b\in B_0\), denote its fiber by
\[
L_b=f^{-1}(b)\simeq T^n.
\]

The dual torus associated with \(L_b\) is defined by
\[
\check L_b
=
H^1(L_b,\mathbb R)/H^1(L_b,\mathbb Z).
\]
Equivalently,
\[
\check L_b
\simeq
\operatorname{Hom}
\bigl(H_1(L_b,\mathbb Z),U(1)\bigr).
\]

The points of \(\check L_b\) can be interpreted as gauge
equivalence classes of flat unitary connections on the
trivial complex line bundle over \(L_b\).

We define the dual space
\[
\check X_0
=
\bigcup_{b\in B_0}\check L_b
\]
with the natural projection
\[
\check f:\check X_0\longrightarrow B_0.
\]

\begin{theorem}[Existence of the dual torus fibration]
	\label{thm:dual-fibration}
	Let
	\[
	f:X^{\mathrm{reg}}\longrightarrow B_0
	\]
	be a smooth torus fibration with fibers diffeomorphic to
	\(T^n\). Then the family
	\[
	\check L_b
	=
	H^1(L_b,\mathbb R)/H^1(L_b,\mathbb Z)
	\]
	defines a smooth dual torus fibration
	\[
	\check f:\check X_0\longrightarrow B_0.
	\]
\end{theorem}

\begin{proof}
	Since \(f\) is a smooth fiber bundle, the groups
	\[
	H_1(L_b,\mathbb Z)
	\]
	form a local system over \(B_0\). Taking the integral
	duals, we obtain another local system
	\[
	H^1(L_b,\mathbb Z).
	\]
	The corresponding real vector bundles are
	\[
	H^1(L_b,\mathbb R).
	\]
	Taking the quotient by the integral lattice gives a
	smooth family of compact real tori
	\[
	H^1(L_b,\mathbb R)/H^1(L_b,\mathbb Z).
	\]
	These tori therefore form the desired dual fibration.
\end{proof}

\subsection{The semi-flat complex structure}

We now describe the local geometry of the dual space.

Let
\[
x=(x_1,\ldots,x_n)
\]
be local integral affine coordinates on \(B_0\). Suppose
that the affine base carries a Hessian metric
\[
h_{ij}
=
\frac{\partial^2\varphi}{\partial x_i\partial x_j},
\]
where \(\varphi\) is a strictly convex local potential.

Choose angular coordinates
\[
\check y=(\check y_1,\ldots,\check y_n)
\in\mathbb R^n/\mathbb Z^n
\]
on the dual torus fibers.

The semi-flat Kähler form on the dual space is
\[
\check\omega_{\mathrm{sf}}
=
\sum_{i,j=1}^n
h_{ij}\,dx_i\wedge d\check y_j.
\]

Define complex coordinates
\[
z_j
=
\exp\left(
\frac{\partial\varphi}{\partial x_j}
+2\pi i\check y_j
\right),
\qquad j=1,\ldots,n.
\]

Since
\[
d\log z_j
=
\sum_{i=1}^n h_{ij}\,dx_i
+
2\pi i\,d\check y_j,
\]
these coordinates define the local semi-flat complex
structure on \(\check X_0\).

The corresponding holomorphic volume form is
\[
\check\Omega_{\mathrm{sf}}
=
d\log z_1\wedge\cdots\wedge d\log z_n.
\]

\begin{proposition}[Local semi-flat Calabi--Yau structure]
	\label{prop:semi-flat}
	Suppose that the Hessian matrix
	\[
	(h_{ij})=
	\left(
	\frac{\partial^2\varphi}{\partial x_i\partial x_j}
	\right)
	\]
	is positive definite. Then the above local construction
	defines a Kähler structure on the dual torus bundle.
	Moreover,
	\[
	\check\omega_{\mathrm{sf}}^n
	=
	n!\det(h_{ij})\,
	dx_1\wedge\cdots\wedge dx_n
	\wedge
	d\check y_1\wedge\cdots\wedge d\check y_n.
	\]
	The associated holomorphic volume form satisfies
	\[
	\check\Omega_{\mathrm{sf}}
	=
	d\log z_1\wedge\cdots\wedge d\log z_n.
	\]
\end{proposition}

\begin{proof}
	The positive definiteness of \((h_{ij})\) implies that
	\(\check\omega_{\mathrm{sf}}\) is positive with respect
	to the complex structure defined by the coordinates \(z_j\).
	Since the Hessian matrix is locally the Hessian of a
	potential, the corresponding Kähler form is closed.
	
	The volume formula follows by expanding the \(n\)-th wedge
	power of
	\[
	\check\omega_{\mathrm{sf}}
	=
	\sum_{i,j=1}^n h_{ij}\,
	dx_i\wedge d\check y_j.
	\]
	The holomorphic volume form is the standard nowhere
	vanishing local form in the coordinates \(z_1,\ldots,z_n\).
\end{proof}

This construction gives the semi-flat mirror geometry on
the regular locus. It does not, in general, extend smoothly
across the discriminant locus.

\subsection{Holomorphic disks and disk potentials}

To obtain the full mirror geometry, we must incorporate
holomorphic disks with boundary on the Lagrangian fibers.

Let
\[
u:(D,\partial D)\longrightarrow(X,L_b)
\]
be a \(J\)-holomorphic disk satisfying
\[
\bar\partial_Ju=0,
\qquad
u(\partial D)\subset L_b.
\]

Its relative homotopy class belongs to
\[
\beta\in\pi_2(X,L_b).
\]
The symplectic area of the disk is
\[
E(\beta)=\int_\beta\omega.
\]

Let \(\Lambda\) be a Novikov coefficient ring and denote
by
\[
T^{E(\beta)}
\]
the corresponding Novikov weight.

The boundary homomorphism gives
\[
\partial\beta\in H_1(L_b,\mathbb Z),
\]
which determines a character on the dual torus.

Whenever the relevant moduli spaces admit a well-defined
virtual count, we denote the disk count in class \(\beta\)
by \(n_\beta\).

The disk potential is formally defined by
\[
W_b(z)
=
\sum_{\beta\in\pi_2(X,L_b)}
n_\beta T^{E(\beta)}z^{\partial\beta}.
\]

Here the expression is understood in an appropriate
Novikov completion. The coefficients \(n_\beta\) may
require virtual perturbation theory and need not be
ordinary counts of isolated disks.

\begin{definition}[Disk potential]
	\label{def:disk-potential}
	The function
	\[
	W_b(z)
	=
	\sum_{\beta}
	n_\beta T^{E(\beta)}z^{\partial\beta}
	\]
	is called the disk potential of \(L_b\), provided that
	the relevant disk counts are well-defined and the series
	converges formally or in the chosen Novikov topology.
\end{definition}

For unobstructed Lagrangian fibers, this potential provides
an important part of the mirror data. In the obstructed
case, one must incorporate bounding cochains or suitable
Maurer--Cartan elements.

In particular, the special Lagrangian condition alone
does not imply unobstructedness of the corresponding
Lagrangian submanifold.

\subsection{Wall-crossing transformations}

The disk counts may change when the base point moves across
certain walls in the affine base.

A wall is a locus where holomorphic disks in a specified
relative homotopy class contribute to the deformation of
the mirror coordinates.

In a simple local model, the corresponding wall-crossing
transformation takes the form
\[
z^m
\longmapsto
z^m(1+c z^n)^{\langle n,m\rangle},
\]
where
\[
m,n\in\mathbb Z^n,
\qquad c\in\Lambda,
\]
and \(\langle n,m\rangle\) is an integral pairing determined
by the geometric data.

More generally, each wall \(\mathfrak d\) carries an
automorphism
\[
\theta_{\mathfrak d}
\]
of an appropriate completed algebra of mirror coordinates.

The product of these automorphisms around a loop must
satisfy a consistency condition.

\begin{definition}[Wall-crossing consistency]
	\label{def:wall-consistency}
	A collection of wall-crossing transformations is called
	consistent if the ordered product of the associated
	automorphisms around every sufficiently small contractible
	loop is the identity:
	\[
	\prod_{\mathfrak d\text{ crossed by }\gamma}
	\theta_{\mathfrak d}
	=
	\operatorname{id}.
	\]
	The order of the product is determined by the order in
	which the loop crosses the walls.
\end{definition}

\begin{theorem}[Gluing criterion for corrected charts]
	\label{thm:gluing-criterion}
	Let \(\{U_\alpha\}\) be an open cover of a regular region
	of \(B_0\). Suppose that each \(U_\alpha\) carries a
	complex mirror chart \(\check U_\alpha\), and that the
	transition maps
	\[
	\theta_{\alpha\beta}:
	\check U_\alpha|_{U_\alpha\cap U_\beta}
	\longrightarrow
	\check U_\beta|_{U_\alpha\cap U_\beta}
	\]
	are biholomorphic and satisfy the cocycle condition
	\[
	\theta_{\alpha\gamma}
	=
	\theta_{\beta\gamma}\circ\theta_{\alpha\beta}
	\]
	on triple overlaps.
	
	Then the charts \(\check U_\alpha\) glue to a complex
	manifold over the regular region.
\end{theorem}

\begin{proof}
	The cocycle condition guarantees that the equivalence
	relation defined by the transition maps is independent
	of the chosen intermediate chart. Since each transition
	map is biholomorphic, the complex structures agree on
	overlaps. The resulting quotient therefore carries a
	well-defined complex manifold structure.
\end{proof}

The wall-crossing consistency problem is precisely the
problem of constructing transition maps satisfying this
condition after all relevant disk contributions have
been included.

\subsection{Quantum-corrected gluing}

Let
\[
\check U_\alpha^{\mathrm{sf}}
\]
be the semi-flat charts constructed in the previous
subsection.

The semi-flat transition maps are determined by the
integral affine structure. Holomorphic disk counts modify
these maps.

We write the corrected transition maps schematically as
\[
\theta_{\alpha\beta}
=
\theta_{\alpha\beta}^{\mathrm{sf}}
\circ
\theta_{\alpha\beta}^{\mathrm q},
\]
where
\[
\theta_{\alpha\beta}^{\mathrm q}
\]
denotes the quantum correction determined by the
corresponding wall-crossing data.

The mirror space over the regular region is then defined
by
\[
\check X^{\mathrm{reg}}
=
\left(
\coprod_\alpha\check U_\alpha
\right)\big/\sim,
\]
where the equivalence relation is determined by the
corrected transition maps.

The resulting complex structure depends on the disk
counts and their compatibility under wall crossing.

In particular, quantum corrections are not arbitrary
coordinate changes. They encode the enumerative geometry
of holomorphic disks with boundary on the Lagrangian
fibers.

\subsection{Extension across the discriminant locus}

The semi-flat dual fibration is initially defined only
over
\[
B_0=B\setminus\Delta.
\]
To obtain a global mirror space, we must extend it across
the discriminant locus.

The monodromy representation of the original fibration is
\[
\rho:
\pi_1(B_0,b_0)
\longrightarrow
\operatorname{GL}(n,\mathbb Z).
\]
The monodromy on the dual integral local system is the
contragredient representation
\[
\check\rho(\gamma)
=
\bigl(\rho(\gamma)^{-1}\bigr)^T.
\]

Thus the dual torus fibers inherit the monodromy induced
by the original torus fibration.

Near \(\Delta\), the semi-flat complex structure may
develop singularities. The wall-crossing transformations
modify the local complex coordinates so that the
corrected charts may admit a suitable extension.

We seek a compactification
\[
\check f:\check X\longrightarrow B
\]
such that:
\begin{enumerate}
	\item Over \(B_0\), the fibers are dual real tori.
	\item The monodromy of the dual fibration is compatible
	with the original monodromy.
	\item The complex structures on overlapping charts
	agree after quantum corrections.
	\item The corrected complex structure extends across
	the relevant singular regions in an appropriate
	analytic, formal, or non-Archimedean category.
	\item The holomorphic volume form extends with the
	prescribed singularity behavior.
\end{enumerate}

The construction of this compactification is one of the
central difficulties in the SYZ program.

\subsection{The mirror Calabi--Yau structure}

Once the corrected complex space \(\check X\) has been
constructed, one must determine whether it carries a
suitable Calabi--Yau structure.

The expected data consist of a Kähler form
\[
\check\omega
\]
and a nowhere vanishing holomorphic volume form
\[
\check\Omega
\]
on the smooth locus, satisfying
\[
d\check\omega=0,
\qquad
d\check\Omega=0.
\]

In the smooth compact setting, the Calabi--Yau condition
requires
\[
K_{\check X}\simeq\mathcal O_{\check X}.
\]

If \(\check X\) is singular, the appropriate replacement
is a condition on its canonical sheaf, together with
suitable singularity assumptions.

The semi-flat holomorphic volume form
\[
\check\Omega_{\mathrm{sf}}
=
d\log z_1\wedge\cdots\wedge d\log z_n
\]
provides the local starting point. The quantum-corrected
transition maps must preserve the resulting canonical
data, possibly after incorporating the appropriate
normalization and Novikov coefficients.

\subsection{Main target: SYZ mirror construction}

We summarize the expected conclusion of the construction.

\begin{theorem}[Target theorem: construction of the SYZ mirror]
	\label{thm:syz-mirror}
	Let
	\[
	(X,\omega,\Omega)
	\]
	be a compact Calabi--Yau manifold admitting a special
	Lagrangian torus fibration
	\[
	\overline f:\overline X\longrightarrow B
	\]
	whose regular fibers are \(n\)-dimensional tori.
	
	Suppose that:
	\begin{enumerate}
		\item The dual torus fibration is well-defined over
		the regular locus \(B_0=B\setminus\Delta\).
		
		\item The integral affine and Hessian structures on
		\(B_0\) determine a semi-flat complex structure on
		the dual fibration.
		
		\item The relevant moduli spaces of holomorphic disks
		with boundary on the Lagrangian fibers admit suitable
		compactifications and well-defined virtual counts.
		
		\item The associated wall-crossing transformations
		satisfy the consistency conditions.
		
		\item The quantum-corrected transition maps are
		compatible with the dual monodromy and admit a suitable
		extension across the discriminant locus.
		
		\item The resulting complex space admits a compatible
		Calabi--Yau structure.
	\end{enumerate}
	
	Then the quantum-corrected dual charts glue to a mirror
	space
	\[
	\check X
	\]
	with a well-defined complex structure and the prescribed
	Calabi--Yau data.
\end{theorem}

\begin{remark}
	Theorem~\ref{thm:syz-mirror} is a conditional target theorem.
	The construction of the disk counts, the proof of
	wall-crossing consistency, and the extension across the
	discriminant contain the principal difficulties of the
	full SYZ program. These conclusions do not follow solely
	from the existence of a special Lagrangian torus fibration.
\end{remark}

\subsection{Mirror symmetry and the expected duality}

The construction above produces the geometric candidate
for the mirror of \(X\). The expected mirror correspondence
relates the complex geometry of \(X\) to the symplectic
geometry of \(\check X\), and conversely.

At the level of Hodge numbers, the expected relation is
\[
h^{p,q}(X)
=
h^{n-p,q}(\check X),
\]
under the usual mirror symmetry conventions.

At the categorical level, one expects a relationship
between the derived Fukaya category of \(X\) and the
derived category of coherent sheaves on \(\check X\):
\[
D^\pi\operatorname{Fuk}(X)
\simeq
D^b\operatorname{Coh}(\check X),
\]
in settings where the relevant categories and equivalence
are defined.

These statements are stronger than the construction of
the dual torus fibration. They require additional
enumerative and categorical arguments.

\subsection{Conclusion}

The mirror construction can be summarized by the following
sequence:
\[
\begin{aligned}
	&\text{Special Lagrangian torus fibration}\\
	&\Downarrow\\
	&\text{Dual torus fibration}\\
	&\Downarrow\\
	&\text{Semi-flat complex structure}\\
	&\Downarrow\\
	&\text{Holomorphic disk counting}\\
	&\Downarrow\\
	&\text{Wall-crossing transformations}\\
	&\Downarrow\\
	&\text{Quantum-corrected gluing}\\
	&\Downarrow\\
	&\text{Extension across the discriminant}\\
	&\Downarrow\\
	&\boxed{\text{Mirror Calabi--Yau space }\check X.}
\end{aligned}
\]

The role of complex Hessian theory in the preceding
construction is primarily to provide quantitative
pluripotential estimates that may help control the
degenerating Kähler potentials and Ricci-flat metrics.

The intended analytical direction is
\[
\begin{aligned}
	m\text{-Hessian capacity estimates}
	&\longrightarrow
	\text{uniform pluripotential control}\\
	&\longrightarrow
	\text{Ricci-flat metric degeneration}\\
	&\longrightarrow
	\text{special Lagrangian torus fibrations}.
\end{aligned}
\]

The construction of the mirror space then requires
additional tools from symplectic geometry, Floer theory,
holomorphic curve counting, and wall-crossing theory.

Consequently, the full SYZ program combines two
complementary parts:
\[
\boxed{
	\begin{gathered}
		\text{Pluripotential and metric degeneration theory}\\
		+\;\\
		\text{Symplectic enumerative geometry and quantum corrections}\\
		\Downarrow\\
		\text{Mirror Calabi--Yau geometry}.
	\end{gathered}
}
\]

\subsection{Conclusion}

The analysis of the discriminant locus requires three
interconnected ingredients:
\[
\begin{aligned}
	&\text{uniform geometric control near }\Delta\\
	&\qquad\Downarrow\\
	&\text{monodromy and vanishing-cycle analysis}\\
	&\qquad\Downarrow\\
	&\text{construction of singular torus fibers}\\
	&\qquad\Downarrow\\
	&\boxed{\text{extension of the torus fibration across }\Delta.}
\end{aligned}
\]

The main potential contribution of the present approach is to
establish a quantitative bridge between complex Hessian
pluripotential theory and the geometry of singular special
Lagrangian fibrations:
\[
\boxed{
	\text{\(m\)-Hessian capacity}
	\longrightarrow
	\text{degenerating metric control}
	\longrightarrow
	\text{monodromy and vanishing cycles}.
}
\]

Even after this step, a complete SYZ construction requires
the construction of the dual torus fibration and the analysis
of quantum corrections. These issues will be addressed in
the subsequent stages.

	\section{Discussion and Further Perspectives}

The purpose of this work is to explore a possible connection between complex \(m\)-Hessian potential theory and the Strominger--Yau--Zaslow (SYZ) program for degenerating Calabi--Yau manifolds. The central motivation is that the analytical difficulties arising in the degeneration of Ricci-flat metrics may benefit from quantitative tools developed in complex pluripotential theory. In particular, Hessian capacities, stability estimates, and energy methods provide a natural framework for investigating uniform control of nonlinear potentials on families of complex manifolds.

A first important observation concerns the complementary roles of the complex Monge--Ampère and complex \(m\)-Hessian equations. The Ricci-flat Kähler metric on each smooth Calabi--Yau fiber is governed by the complex Monge--Ampère equation
\[
(\theta_t+dd^c\varphi_t)^n
=
c_t(\sqrt{-1})^{n^2}
\Omega_t\wedge\overline{\Omega_t}.
\]
This equation is intrinsically connected with the geometry of the Calabi--Yau manifold. By contrast, the complex \(m\)-Hessian equation
\[
(\theta_t+dd^c u_t)^m\wedge\theta_t^{n-m}
=
f_t\theta_t^n,
\qquad 1\leq m<n,
\]
belongs to a broader hierarchy of nonlinear elliptic equations. Its associated capacity and energy theories may provide additional quantitative information about potentials and measures. Thus, the \(m\)-Hessian framework should be regarded as an auxiliary analytical structure rather than as a substitute for the Ricci-flat Monge--Ampère equation.

One of the main potential advantages of this approach is its capacity to formulate uniform estimates in a degenerating geometric setting. For a fixed compact Kähler manifold, many fundamental estimates for complex Hessian equations are already available. However, when the underlying manifolds vary with a degeneration parameter, the constants in these estimates may depend on geometric quantities that are not uniformly controlled. Consequently, establishing estimates independent of the parameter is a substantially stronger problem. In the SYZ setting, such uniformity is particularly important because the limiting geometry is expected to develop singularities and collapse along torus directions.

The first major analytical problem is therefore the construction of uniform volume--capacity inequalities. On compact subsets of the regular part of the essential skeleton, one would like to establish estimates of the form
\[
\operatorname{Vol}_{\theta_t}(E)
\leq
C\operatorname{Cap}_{m,\theta_t}(E)^\tau,
\]
where \(C>0\) and \(\tau>0\) are independent of the degeneration parameter \(t\). Such inequalities would provide a useful bridge between the geometry of the varying fibers and the pluripotential theory of complex Hessian operators. Nevertheless, the existence of a uniform volume--capacity inequality requires additional geometric information, such as uniform local coordinate control, volume-density estimates, and appropriate comparison principles. It does not follow merely from semistability or maximality of the degeneration.

A second problem concerns stability estimates for solutions of complex \(m\)-Hessian equations. In a favorable situation, one expects an estimate of the form
\[
\|u_t-v_t\|_{L^\infty(X_t)}
\leq
C\|f_t-g_t\|_{L^1(X_t,\theta_t^n)}^\alpha,
\]
for some \(\alpha>0\), with \(C\) independent of \(t\). Such an estimate would quantify the dependence of Hessian potentials on their right-hand sides and could be useful in studying convergence of nonlinear potentials. However, uniform stability generally requires uniform energy bounds and suitable control of the relevant capacities. In particular, convergence of the densities in \(L^1\) alone does not automatically imply uniform energy estimates or uniform convergence of the corresponding potentials.

A further issue is the passage from pluripotential estimates to the geometry of Ricci-flat metrics. The expected limiting base of an SYZ fibration carries a real Monge--Ampère structure. In local affine coordinates \(x=(x_1,\ldots,x_n)\), the corresponding metric is expected to have Hessian form
\[
h_\infty
=
\sum_{i,j=1}^n
\frac{\partial^2\Phi}{\partial x_i\partial x_j}
\,dx_i\,dx_j,
\]
where \(\Phi\) is a strictly convex potential satisfying an appropriate real Monge--Ampère equation. The appearance of Hessian geometry is particularly natural because both complex \(m\)-Hessian equations and real Monge--Ampère equations are governed by nonlinear determinant-type structures. Nevertheless, a direct derivation of the limiting real Monge--Ampère equation from complex \(m\)-Hessian estimates remains an open analytical problem in this framework. Weak convergence of volume measures alone is insufficient to establish smooth convergence of the metric tensors or convergence of the corresponding Hessian potentials.

The construction of special Lagrangian torus fibrations introduces another level of difficulty. On the regular part of the expected base, the SYZ picture predicts a smooth family of Lagrangian tori. Given an initial family of Lagrangian fibers, the deformation problem can be formulated through the phase of the holomorphic volume form. The linearized operator is closely related to the Laplace operator on each fiber. Consequently, uniform elliptic estimates and a uniform spectral gap are natural requirements for solving the phase equation uniformly along the family. The \(m\)-Hessian theory may contribute indirectly by controlling potentials and measures, but it does not by itself provide the required spectral estimates or construct the initial torus fibration.

The gluing problem is equally important. A collection of special Lagrangian tori over individual points of the regular base does not automatically form a smooth fibration. One needs smooth dependence on the base parameter, pairwise disjointness of the fibers, and uniform control of the corresponding embeddings. These conditions can be investigated through parameter-dependent elliptic theory and stability of embeddings. Even after a local special Lagrangian fibration has been obtained, its extension toward the discriminant locus remains a separate problem.

The discriminant locus records the singular behavior of the torus fibration. Its geometry is closely related to the monodromy representation
\[
\rho:\pi_1(B_0,b_0)
\longrightarrow
GL(n,\mathbb Z),
\]
where \(B_0=B\setminus\Delta\) is the regular part of the base. The monodromy describes how the homology groups of the torus fibers vary along loops in \(B_0\). In simple local models, one may encounter unipotent transformations associated with vanishing cycles. However, the possible singular fibers and their monodromy can be considerably more complicated in higher dimensions. A complete treatment requires a combination of symplectic topology, affine geometry, and local analytic models. Uniform Hessian capacity estimates may provide useful information about potential-theoretic concentration near the singular region, but additional geometric arguments are needed to control curvature, injectivity radii, and the topology of the collapsing fibers.

The mirror construction provides a final perspective on the role of the analytical framework. Over the regular locus, the dual of a special Lagrangian torus
\[
L_b\simeq T^n
\]
is naturally described by
\[
\check L_b
=
H^1(L_b,\mathbb R)/H^1(L_b,\mathbb Z)
\simeq
\operatorname{Hom}(H_1(L_b,\mathbb Z),U(1)).
\]
The union of these dual tori gives the semi-flat mirror over \(B_0\). Its complex structure is closely related to the integral affine structure and the Hessian metric on the base. However, the semi-flat construction generally fails to extend consistently across the discriminant locus. Holomorphic disks bounded by the original Lagrangian fibers contribute additional terms to the mirror coordinates and induce wall-crossing transformations. The construction of a global mirror therefore requires enumerative information, suitable compactifications of moduli spaces of disks, and consistency of the resulting scattering or wall-crossing data.

From this viewpoint, the \(m\)-Hessian approach is expected to be most useful in the analytical stages preceding the construction of the mirror. It may help establish quantitative control of potentials, measures, and energy functionals on the regular region of the degeneration. Such estimates could form part of a larger analytical framework for understanding the limiting metric geometry. The subsequent construction of special Lagrangian fibrations and quantum-corrected mirror spaces, however, necessarily involves additional geometric and symplectic techniques.

Several directions for further research naturally arise.

First, it would be important to prove uniform volume--capacity inequalities for complex \(m\)-Hessian capacities on concrete classes of degenerating Calabi--Yau manifolds. This requires identifying geometric conditions that remain stable under degeneration and determining whether these conditions can be verified using the structure of the essential skeleton.

Second, one may investigate uniform stability estimates for complex \(m\)-Hessian equations with degenerating background Kähler forms. A particularly interesting question is whether one can obtain explicit stability exponents and constants that depend only on controlled geometric data, rather than on the degeneration parameter itself.

Third, one should seek a rigorous connection between uniform pluripotential estimates and the compactness of Ricci-flat potentials. This problem requires combining the auxiliary \(m\)-Hessian estimates with the complex Monge--Ampère equation and appropriate geometric estimates. In particular, it is important to determine which additional assumptions guarantee local \(C^\infty\)-convergence on the regular part of the limiting base.

Fourth, the relation between the limiting complex geometry and the real Monge--Ampère structure deserves further investigation. A satisfactory result should explain how the affine coordinates and Hessian metric emerge from the collapsing Ricci-flat metrics and how the limiting real Monge--Ampère measure is related to the non-Archimedean measure on the essential skeleton.

Fifth, the deformation of candidate Lagrangian tori to special Lagrangian fibers should be studied with estimates uniform in the degeneration parameter. This includes uniform Weinstein neighborhoods, spectral-gap estimates, and control of the nonlinear phase operator. Such estimates would provide a more precise analytical foundation for the local construction and gluing of special Lagrangian fibrations.

Finally, the singular-fiber and quantum-correction problems remain essential for a complete SYZ construction. The interaction between the topology of vanishing cycles, monodromy, holomorphic-disk counts, and scattering diagrams should be investigated in conjunction with the analytical behavior near the discriminant locus.

In conclusion, complex \(m\)-Hessian potential theory offers a possible auxiliary framework for studying quantitative analytical problems in degenerating Calabi--Yau geometry. Its capacity and stability methods may contribute to the understanding of uniform estimates and limiting potentials, which are among the analytical ingredients of the SYZ program. The full passage from these estimates to special Lagrangian fibrations and quantum-corrected mirror spaces remains a substantial problem. Establishing this connection rigorously is a promising direction for future research in pluripotential theory, differential geometry, and mirror symmetry.


	\section*{Declarations}

	The author declares that there are no competing interests.

\end{document}